\documentclass[12pt,a5paper]{article}

\pdfoutput=1
\usepackage[pdftex]{hyperref}
\usepackage{color,pdfcomment}
\usepackage{euler,amsmath,amsfonts,url,doi,graphicx}
\usepackage{jeb}

\def\sgn{\operatorname{sgn}}

\usepackage{pgfplots,tikz} 
\usepgfplotslibrary{external} 
\usetikzlibrary{decorations.markings}
\usetikzlibrary{arrows}

\allowdisplaybreaks

\author{J.~E. Bunder\thanks{School of Mathematical Sciences, University of Adelaide, South Australia~5005, Australia.  
\protect\url{mailto:judith.bunder@adelaide.edu.au}} 
\and A.~J. Roberts\thanks{School of Mathematical Sciences, University of Adelaide, South Australia~5005, Australia.  
\protect\url{mailto:anthony.roberts@adelaide.edu.au}}
\and I.~G.~Kevrekidis\thanks{Department of Chemical and Biological Engineering and PACM, Princeton University, Princeton, NJ~08544, USA.}}

\date{\today}

\title{Accuracy of patch dynamics with mesoscale temporal coupling for efficient exascale simulation}

\begin{document}
\maketitle

\begin{abstract}
Massive parallelisation has lead to a dramatic increase in available computational power. 
However, data transfer speeds have failed to keep pace and are the major limiting factor in the development of exascale computing.
New algorithms must be developed which minimise the transfer of data.
Patch dynamics is a computational macroscale modelling scheme which provides a coarse macroscale solution of a problem defined on a fine microscale by dividing the domain into many nonoverlapping, coupled patches. 
Patch dynamics is readily adaptable to massive parallelisation as each processor can evaluate the dynamics on one, or a few, patches.  
However, patch coupling conditions interpolate across the unevaluated parts of the domain between patches, and are typically reevaluated at every microscale time step, thus requiring almost continuous data transfer.
We propose a modified patch dynamics scheme which minimises data transfer by only reevaluating the patch coupling conditions at `mesoscale' time scales which are significantly larger than the microscale time of the microscale problem.
We analyse the error arising from patch dynamics with mesoscale temporal coupling as a function of the mesoscale time interval, patch size, and ratio between the microscale and macroscale. 
\end{abstract}

\tableofcontents

\section{Introduction}
\label{sec:intro}

Mathematical equations describing a phenomenon (e.g., fluid flow, chemotaxis, mechanics) are typically written at the scale at which we want the information (e.g., macroscopic velocity fields, bacterial concentrations, macroscopic deformations). 
Increasingly, the scale at which the physics are understood (molecular, cellular, agent-based) is much finer than the macroscopic, human, systems scale at which we want information and a variety of modelling techniques may be applied to reinterpret the microscale problem at the desired macroscale~\cite[e.g.]{Attinger2004,Givon04, Horstemeyer2010, Dada2011, Degenhard2006}. 
Unfortunately, for many multiscale and multiphysics problems good macroscale closures do not exist---instead we simulate and observe very detailed models of great complexity at great cost \cite[e.g.]{Ortega2013,Nguyen2014,Kiuchi2014,Miyoshi2014}. 
In this scenario, we aim to develop efficient computational procedures to wrap around whatever microscopic level computer model a scientist chooses for any given system \cite[e.g.]{Kevrekidis03b, Kevrekidis04a}---be it anything from a Monte--Carlo description of a chemical reaction to an individual\slash agent-based model in ecology or epidemiology.

The methodology is to evaluate automatically (`on demand'), directly from the micoscale model, the macroscopic modelling closures for the emergent dynamics which all too often are not available explicitly \cite[e.g.]{Cisternas03}.
This is an `equation-free' method in the sense that it makes no attempt to derive a macroscale equation, in contrast to, for example, homogenization~\cite{Pavliotis2008}. 

\begin{figure}
\centering
\includegraphics{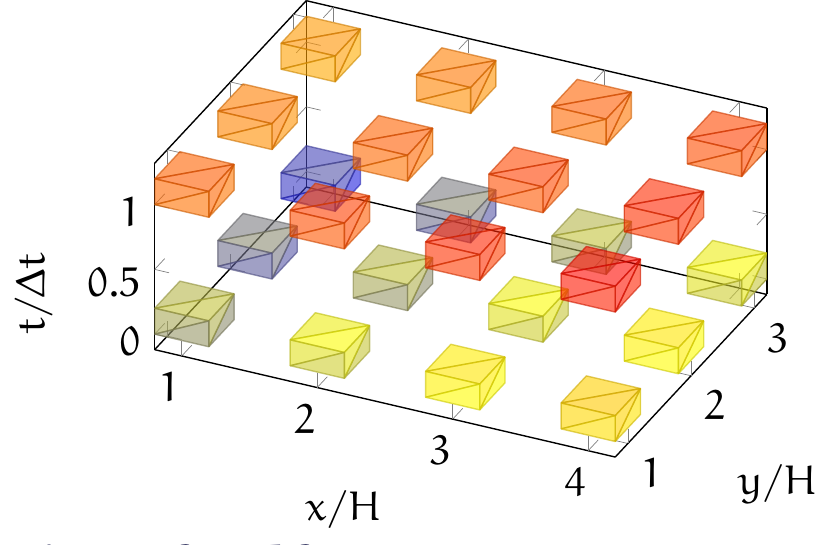}
\caption{Equation-free methods simulate only on small patches in space-time: patches are placed at macroscale time steps~$\Delta t$ and spatial macroscale steps~$H$.  
The given microscale dynamics are only simulated within each patch: coupling conditions interpolate across un-simulated space; and projective integration steps across time.}
\label{fig:patchesSpaceTime}
\end{figure}

The aim of this article is to develop and support the patch dynamics scheme~\cite[for reviews]{Givon04, Hyman05, Li2007, Kevrekidis09, Samaey10}.
By only computing on a small fraction of the space-time domain, see the schematic Figure~\ref{fig:patchesSpaceTime}, this scheme empowers large scale simulation and prediction.
But special challenges arise in the largest simulations on exa\slash peta-scale computers.
Designed for exa\slash peta-scale computing we propose new infrequent couplings between these microscale simulations on microscale patches across un-simulated space.
We establish new results on efficiency, accuracy, and consistency for the emergent macroscale simulation.
Section~\ref{sec:patchsetup} discusses the mathematical details of patch dynamics with patches defined in space only (also referred to as the gap-tooth method~\cite{Roberts07}) and uses the example of a one dimensional diffusion problem on a discrete lattice.

\begin{figure}
\centering
\includegraphics{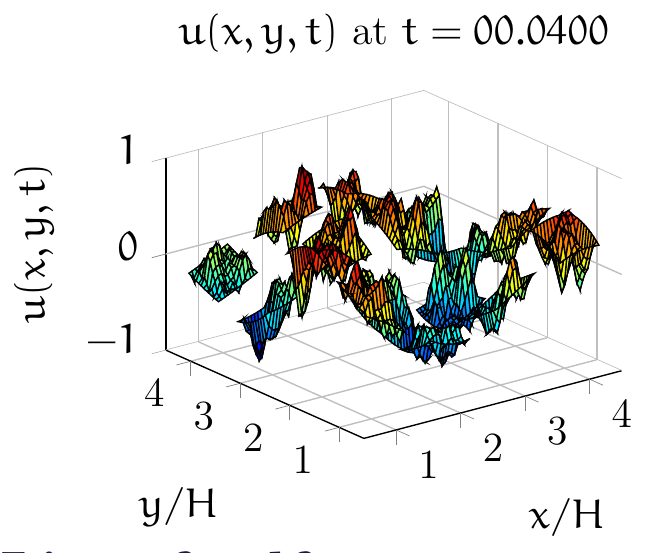}\includegraphics{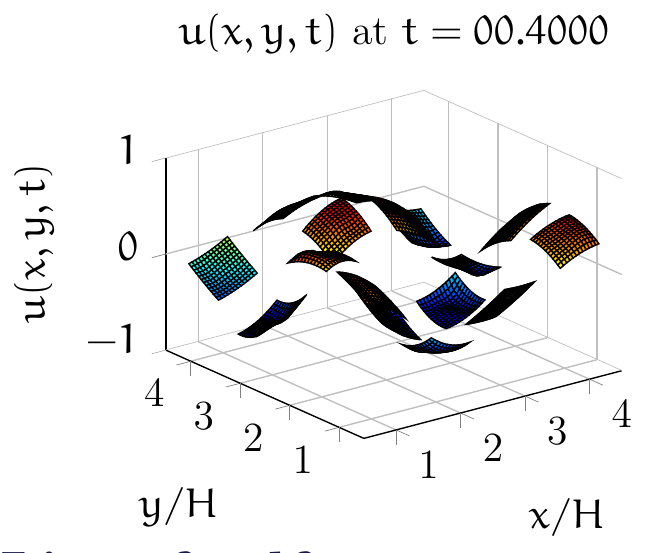}
\includegraphics{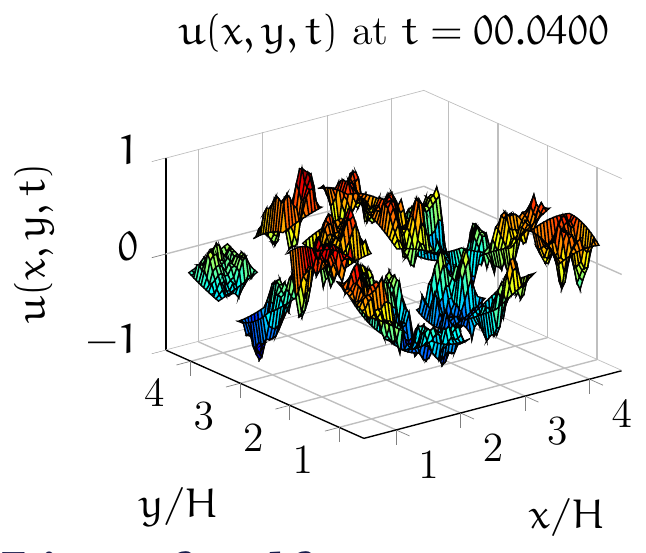}\includegraphics{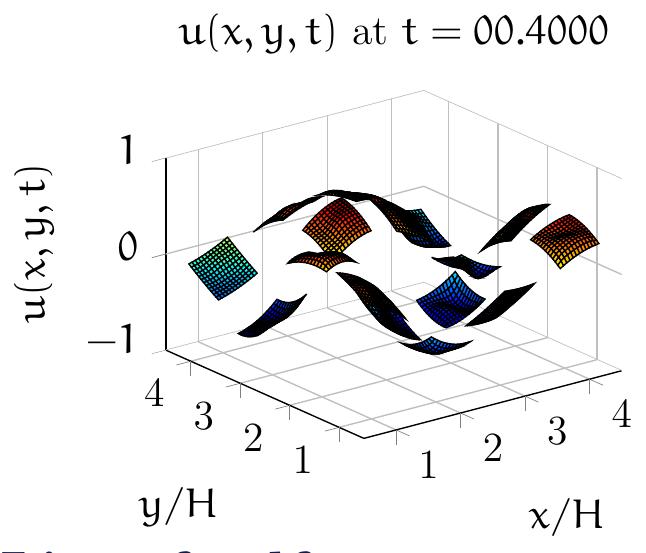}
\caption{Real parts of Ginzburg--Landau microscale fields~$u_{j_x+i_xN,j_y+i_yN}$ in patches with $n=6$ at time~$t=0.04$ and $t=0.4$ with: (top)~continuous time coupling; and (bottom)~infrequent mesoscale coupling $\delta t=0.2$\,.
At this scale there is little to distinguish the continuous time coupling solution and the mesoscale coupling solution.
}\label{fig:3d}
\end{figure}

The scheme is readily adaptable to higher dimensions and more complex nonlinear systems. 
Section~\ref{sec:numerics} develops numerical simulations of a two dimensional complex Ginzburg--Landau \textsc{pde}.
Simulations, such as that shown in Figure~\ref{fig:3d}, qualitatively confirm the accuracy of our proposed infrequent coupling scheme for this two dimensional nonlinear system.

As computational power approaches the exascale,  roughly $10^{18}$\,\textsc{flops} (floating point operations per second), it is tempting to think that soon many multiscale problems will be solved numerically by computing a microscale simulation across the entire domain.
However, constraints on high performance computing make such a task effectively impossible for all but simple scenarios~\cite{dolbow2004}.
Improvements in high performance computing are gained through massive parallelisation via increases in the number of  processors, but  success is forecast to be limited~\cite{Kogge2008, doe2009, Ashby2010, doe2013, emwg2014}.
First, one limitation is that memory storage is growing at a tenth of the pace of processing power~\cite{Ashby2010}, so there is huge processing power but little space to store the resulting data.
Schemes for exascale computing must use only as much data as required; that is, we need as sparse as possible a resolution over the macroscale, such as that offered by the patch scheme.
Secondly, relative to computational speeds, the slow speed of data transfer between processors, cache and memory prevents processors from operating effectively unless the computational scheme limits communication, as we propose and analyse here for the patch scheme.
Thirdly, with millions of processors, hardware failure will be common somewhere and we discuss possibilities for fault tolerance in the patch scheme.
These limitations are not expected to be overcome through improvements in hardware, so it is mainly through the development of new algorithms that engineers and scientists will exploit the benefits of massive parallelisation~\cite{doe2010, Ashby2010, doe2013, emwg2014}.

\begin{figure}
\centering
\includegraphics{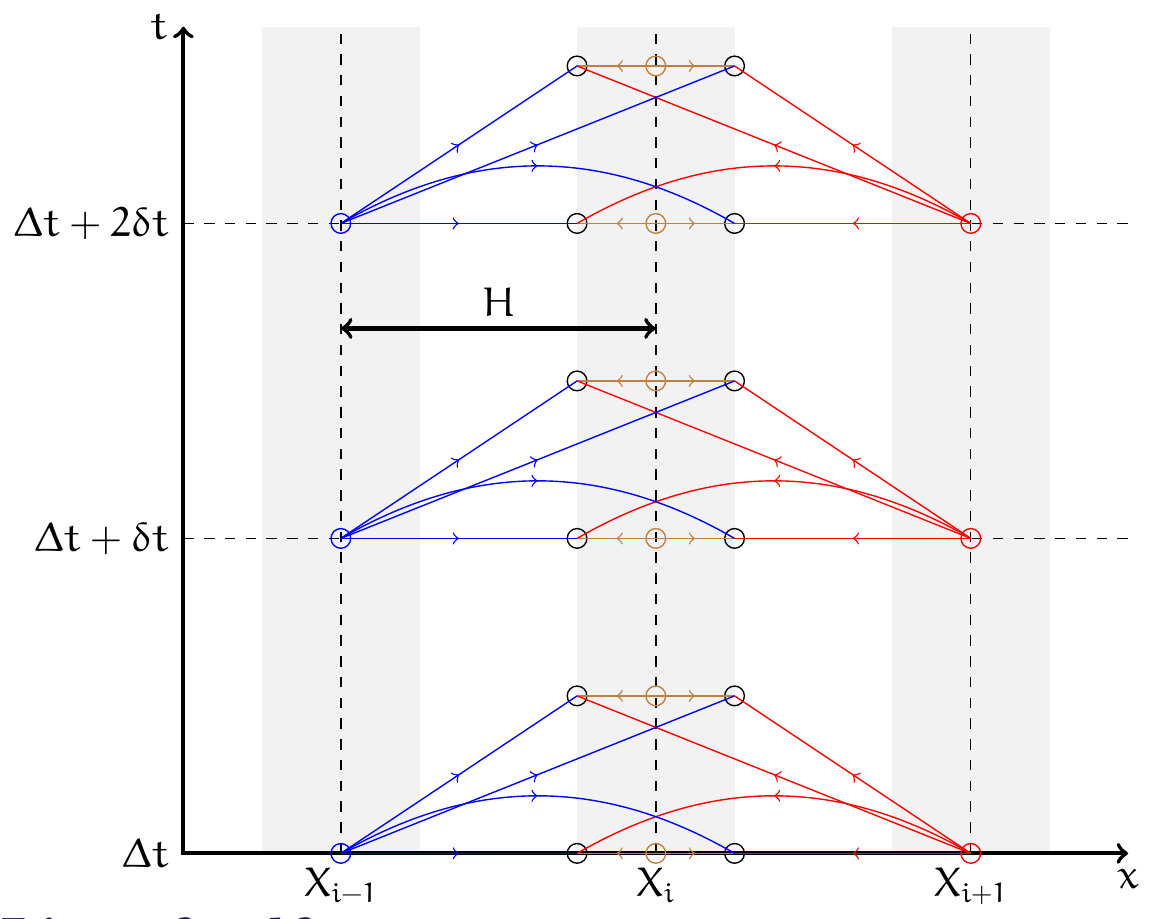}
\caption{Schematic description of the coupling conditions \eqref{eq:cc2}~or~\eqref{eq:cc3} for the $i$th~patch with nearest neighbour coupling where coupling between patches is only reevaluated at mesoscale time steps~$\delta t$.
Three patches are indicated by shaded regions centred about macroscale lattice points~$X_{i-1}$\,, $X_i$ and $X_{i+1}$\,.
The macroscale lattice spacing is~$H$.
As indicated by the coloured lines, the average dynamics on patches $i\pm 1$~and~$i$~\eqref{eq:amp} feed into the coupling conditions \eqref{eq:cc2}~or~\eqref{eq:cc3} on the $i$th~patch and this controls the dynamics on both edges of the $i$th~patch. 
At time~$m\delta t +t$, for nonnegative integer~$m$ and $\Delta t\leq t<\Delta t+\delta t$\,, the coupling conditions of the $i$th~patch are dependent on the average dynamics of the $i$th~patch at time~$m\delta t+t$ (brown) as well as the average dynamics of the neighbouring patches from the mesoscale time step~$m\delta t$ (blue, red).}
\label{fig:elements}
\end{figure}

The patch dynamics approach does not invoke a macroscale equation and requires no prior analysis of the spatial-temporal domain~\cite[e.g.]{Kevrekidis09}. 
This removes significant data storage constraints while increasing the flexibility, allowing `on-the-fly' modifications.
The discretisation of the domain into patches makes patch dynamics readily adaptable to massive parallelisation: for example, a domain decomposition where each processor simulates the dynamics on a few patches.
Here, for simplicity, we assume only one patch per processor.
However, extant implementations of patch dynamics require that coupling conditions are calculated and communicated between each patch at each microscale time step~\cite{Hyman05, Roberts01, Roberts03, Roberts07, Roberts11}, thus requiring substantial and effectively prohibitive data transfer between processors.
Section~\ref{sec:patchsetup} proposes a new modification to the patch scheme to reduce data transfers for exascale computing by limiting the times at which the inter-patch coupling conditions are updated.
As illustrated by the red and blue lines in Figure~\ref{fig:elements}, coupling condition data required from other patches (or processors) is updated only at mesoscale time steps~$\delta t$ which are significantly larger than the microscale time steps of the simulator but smaller than the macroscale time of interest~$\Delta t$. 
However, as indicated by the brown line in Figure~\ref{fig:elements}, the data required for one patch's coupling conditions which is dependent on the dynamics of that patch is updated at microscale time intervals since this information is readily available to the processor.  
These mesoscale coupling adjustments to the patch dynamics scheme should greatly increase the speed of a simulation run on a high performance computer.

A significant issue for exascale computing  is fault management and resiliency~\cite{doe2010, doe2013, emwg2014}.
Typically, if a computer component fails and data is lost, then, assuming the computer is still operational, the required calculation is redone either from scratch or, if there is some fault tolerance written into the algorithm, from a checkpoint.
If a failure causes data to be delayed rather than lost, then the whole computation is delayed until the data is successfully transferred. 
In either case, failures increase the time required to complete a calculation.
Such delays are not usually major issues for a computer with relatively few components, but on an exascale computer with millions of processors and numerous other components, failure is expected to be a regular occurrence and restarting from a checkpoint and waiting for data is not viable~\cite{Schroeder2007, Ashby2010, emwg2014}. 
Algorithms for exascale computing must be fault tolerant while also accounting for errors associated with failure. 
Section~\ref{sec:conclude} discusses how fault management may be incorporated into the proposed patch dynamics scheme.

To demonstrate how to apply patch dynamics mesoscale coupling, Section~\ref{sec:micro} solves a fundamental microscale discrete diffusion problem using standard patch dynamics macroscale modelling without mesoscale coupling; that is, with inter-patch coupling at microscale time intervals. 
Then, Section~\ref{sec:meso} modifies the solution to allow for patch coupling at only mesoscale time intervals.
Section~\ref{sec:error} analyses the error of the solution with mesoscale coupling obtained in Section~\ref{sec:meso} relative to the solution obtained without mesoscale coupling  obtained in Section~\ref{sec:micro}.
Therefore, this error is not the full error of the modelling but requires additional consideration of the error associated with the standard patch dynamics scheme. 
The error of standard patch dynamics has been discussed for several different systems with a variety of microscale structures~\cite{Roberts01, Roberts03, Samaey05, Bunder2012, Roberts11, codevariD}.
Although we only consider the mathematical details for one patch, it is readily scalable to a multiple processor system, as shown in the numerical results of Section~\ref{sec:numerics}.

\section{Patch dynamics implementation}
\label{sec:patchsetup}

As an initial prototype problem we consider a simple diffusion system on a discrete one dimensional microscale spatial lattice with lattice index~$j$ and microscale lattice spacing~$h$, and time~$t$ which is measured on some microscale,
\begin{equation}
\dot{u}_j(t)=u_{j+1}(t)+u_{j-1}(t)-2u_j(t)\,,\label{eq:ode}
\end{equation}
with some given initial condition~$u_j(0)$ on the microscale field. 
Realistic microscale dynamics are much more complicated than this simple microscale diffusion but before we can contemplate realistic dynamics we prove that the proposed procedure for exascale computing is sound for at least this foundational case of the system~\eqref{eq:ode}.
Section~\ref{sec:numerics} presents successful numerical simulations of the more complex two dimensional Ginzburg--Landau \textsc{ode}. 

\begin{figure}
\centering
\includegraphics{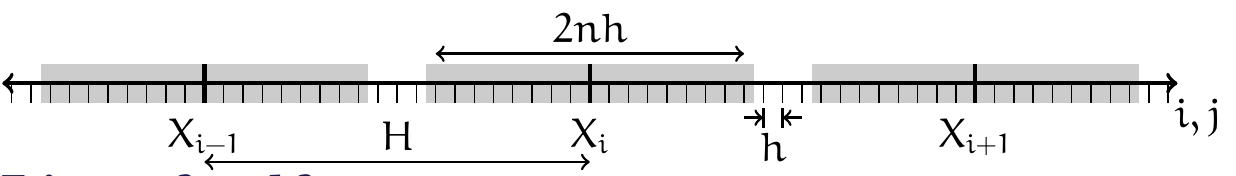}
\caption{The microscale lattice with microscale lattice spacing~$h$, the macroscale lattice with $N=20$ microscale lattice points in one macroscale step, and three patches with patch half-width $n=8$\,.  
The relatively large ratio $n/N=0.4$ means there is little space between the patches.}
\label{fig:patches}
\end{figure}

Suppose we require a macroscale simulation of \textsc{ode}~\eqref{eq:ode} but only at discrete spacings $H\gg h$\,.
We construct a macroscale lattice with spacing~$H$ and macroscale lattice sites~$X_i=iH=iNh$, for patch index $i=0,\pm1,\ldots$ and where the number of microscale lattice points within one macroscale step $N=H/h$ is integral.
In the patch dynamics scheme, for all~$i$ we construct the $i$th~patch of width~$2nh$, for positive integer~$n$, centred about the macroscale lattice site~$X_i$\,:  $n/N<1/2$  ensures the patches do not overlap;  in practice, $n/N\ll 1/2$ for efficient macroscale simulation. 
In Figure~\ref{fig:elements} the shaded areas schematically represent three patches at $X_{i\pm 1}$ and $X_i$\,, and in Figure~\ref{fig:patches} these same patches are superimposed on both the microscale and the macroscale lattices. 
Integer $n$~is the patch half-width; that is, \(n\)~is the number of microscale lattice points which fit into half a patch. 
Figure~\ref{fig:patches} shows patches with patch half-width $n=8$ and $n/N=0.4$\,.
The ratio~$n/N$ is equal to the ratio of half the physical patch width~$nh$ to the macroscale lattice spacing~$H$.

\begin{figure}
\centering
\includegraphics{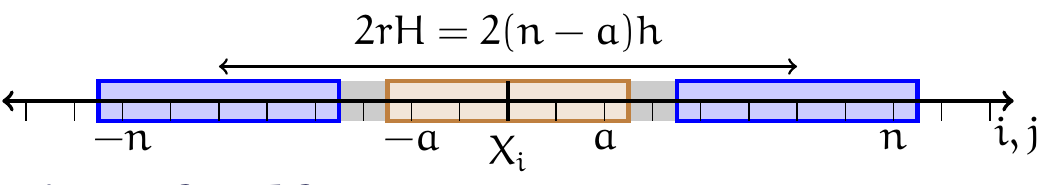}
\caption{A closeup example of the $i$th patch, similar to that shown in Figure~\ref{fig:patches} with patch half-width $n=8$\,.
The core and action half-width is $a=2$\,, and the averaging over the core (brown) provides the macroscale field~\eqref{eq:amp}.
The left and right action regions (blue) for the coupling conditions~\eqref{eq:cc}, \eqref{eq:cc2}~or~\eqref{eq:cc3} are averages over $j=\pm(8,7,\ldots,4)$.
The ratio $r=(n-a)/N=0.3$\,, using $N=20$ from Figure~\ref{fig:patches}.
The buffer width is $n-a=6=rN$\,.
While it is only the core which contributes to the evaluation of the macroscale field~$U_i(t)$, extending the domain of the patch into buffer regions improves the accuracy of the macroscale field. 
}
\label{fig:onepatch}
\end{figure}

We solve \textsc{ode}~\eqref{eq:ode} for microscale fields~$u_j(t)$ but only for the microscale lattice points which lie within a patch; that is, for~$u_{j+iN}$ with microscale sub-patch index $j=0,\pm1,\ldots,\pm (n-1)$ and macroscale patch index $i$.
Without loss of generality we start at $t=0$\,, although the results presented here apply to any initial time which is an integer multiple of mesoscale time step~$\delta t$.
We assume the initial and final simulation times are in the same patch, such as those times shown in Figure~\ref{fig:patches}, so do not consider patch coupling across time. 
The patch boundary conditions of the microscale simulators are provided by patch coupling conditions which extrapolate across the un-simulated space between patches.
Coupling between adjacent patches is achieved by constraining the average of the microscale field near both edges of each patch, termed the `action regions'. 
The left and right action regions of one patch are shaded blue in Figure~\ref{fig:onepatch} and extend over microscale sub-patch indices $j=\pm n,\ldots,(n-2a)$ for some integer $a$, $0\leq a<n$\,.
After solving for all microscale fields within the $i$th~patch we extract the desired macroscale solution~$U_i(t)$ at~$X_i$ via some averaging over the microscale solutions in the middle or `core'  of the $i$th~patch.
The core of one patch is shaded brown in Figure~\ref{fig:onepatch} and extends over $j=0,\pm1,\ldots,\pm a$\,.
The integer~$a$ is defined as the core half-width.
Li et~al.~\cite{Li1999} used similar averaging techniques over action regions and core to  simulate molecular dynamics in a fluid.

The fields~$u_{j+iN}$ in the core (i.e., $j=0,\pm1,\ldots,\pm a$) of each patch~$i$  are the most useful as they define the macroscale solution~$U_i(t)$ which is used for both interpolation between patches and extrapolation across time. 
In contrast, the fields~$u_{j+iN}$ which lie inside a patch but outside the core (i.e., $j=\pm(a+1),\ldots,\pm n$) are of no interest and are forgotten as soon as they are calculated.
However, this region inside a patch but outside the core performs an important function in that it `sheilds' or `buffers' the fields within the core from errors which arise on the application of coupling conditions on the microscale solutions in the action regions.  
The left and right `buffers' extend from the patch edge to the core edge over sub-patch indices $j=\pm(a+1),\ldots,\pm n$ and have width~$n-a$~\cite{codevariD}.
Generally a larger buffer results in smaller numerical error; Figure~\ref{fig:er} in Section~\ref{sec:error} shows errors decreasing approximately exponentially with buffer width.
Thus the buffers perform an important but ancillary task, improving the macroscale solution without directly contributing to its evaluation. 
Section~\ref{sec:micro} shows that $n-a$, rather than~$n$, determines eigenmodes of the microscale solution in one patch and so, in addition to being the buffer width, $n-a$ plays the role of an effective patch half-width.
We name~$n-a$ the reduced patch half-width.

Patch coupling conditions and the method of deriving the macroscale solution may vary, depending on the problem being considered.
In choosing these rules the form of the microscale model is important. 
For example, Section~\ref{sec:error} shows that numerical error decreases with increasing reduced patch half-width~$n-a$, so it is recommended that~$n-a$ is as large as possible.
This recommendation is suitable for many smooth microscale models such as the simple diffusion problem~\eqref{eq:ode} and the two dimensional Ginzburg--Landau equation discussed in Section~\ref{sec:numerics}.
However,  when the microscale model has some rough fine scale structure, such as a periodic spatial roughness, more care needs to be taken when choosing~$n-a$. 
In the case of some periodic spatial roughness, minimal errors are obtained when patch dynamics adequately averages over the spatial structures; for example, when the period exactly divides the reduced patch half-width~$n-a$~\cite{codevariD}.
Since \textsc{ode}~\eqref{eq:ode} has no rough microscale structure we do not need to consider this complication.

The macroscale field obtained from each patch is generally some average over the microscale fields in the centre core of the patch.
For core half-width~$a$, $0\leq a<n$\,, we average over the $2a+1$~microscale fields in the core of the $i$th~patch to define the macroscale field of that patch as
\begin{equation}
U_i(t)=\sum_{j=-a}^a\frac{u_{j+iN}(t)}{2a+1}\,.\label{eq:amp}
\end{equation}
Figure~\ref{fig:onepatch} shows one patch with a core half-width $a=2$\,.
In two or more dimensions the core may be a rectangle or rectangular prism in the centre of a patch and, as in one dimension, the macroscale field of a patch is the average of those microscale fields within the core. 
Section~\ref{sec:numerics} implements square patches with square cores.

Each patch is coupled to its near neighbours via a control acting on the the boundary action regions of each patch.
We define left and right action regions on the patch edges, both containing $2a+1$~microscale lattice points, with $j=-n,-n+1,\ldots,-(n-2a)$ and $j=n,n-1,\ldots,(n-2a)$ respectively.  
With coupling at microscale times, the coupling conditions of the $i$th~patch are proposed to be that the action region average~\cite{codevariD}
\begin{equation} 
\sum_{j=\pm n-2a}^{\pm n}\frac{u_{j+iN}(t)}{2a+1}
=U_i(t)\cos\ell +\frac{f_{i,\pm n}(t)}{2a+1}\,,\label{eq:cc}
\end{equation}
for some $f_{i,\pm n}$ and $\cos\ell$ derived via classic Lagrange interpolation of neighbouring patch macroscale fields~\(U_{i\pm 1}(t)\), \(U_{i\pm 2}(t)\),\ldots and the patch macroscale field~\(U_{i}(t)\) (this interpolation has been proved to be effective for \textsc{pde}{}s~\cite{Roberts07}).
The right hand side of the coupling conditions~\eqref{eq:cc} only contains macroscale fields obtained from applying averaging~\eqref{eq:amp} to microscale fields within the core.
Thus, as discussed above, microscale fields within the buffers make no direct contribution to the interpolation across the un-simulated space between patches.
The action regions are the same size as the core and $a$~is both the core half-width and the action half-width, as illustrated in Figure~\ref{fig:onepatch}.


The details of the patch coupling are contained in $\cos\ell$~and~$f_{i,\pm n}$  which are functions of a parameter~$\gamma$, which controls the coupling strength between patches, and the ratio $r=(n-a)/N$\,, which compares half the physical reduced patch width $(n-a)h$ with the macroscale lattice spacing~$H$. 
For our purposes, the details of the coupling contained in $\cos\ell$~and~$f_{i,\pm n}$ are not important.
Section~\ref{sec:error} shows that errors associated with mesoscale coupling are not dependent on~$\cos\ell$ and the
 error analysis is presented in terms of units of $f_{i,\pm n}$ and its temporal derivatives. 
The derivation of $\cos\ell$~and~$f_{i,\pm n}$ to any order of coupling (i.e., the number of patches coupled to any one patch) is presented elsewhere~\cite{Roberts11, codevariD}.
For example, for only nearest neighbour coupling these function have a parabolic dependence on the ratio~$r$ and linear dependence on coupling parameter $\gamma$: $\cos\ell= (1-r^2\gamma)$ and
\begin{equation}
 f_{i,\pm n}(t)=\tfrac12(2a+1)r\gamma[(r\pm 1)U_{i+1}(t)+(r\mp 1)U_{i-1}(t)]\,,\label{eq:nn}
\end{equation}
where the scaling by~$2a+1$ and the subscript~$\pm n$ on $ f_{i,\pm n}$ are for later convenience.
Since $0<r<1/2$ and $0\leq\gamma\leq 1$\,, for nearest neighbour coupling $0.75<\cos\ell\leq 1$\,.
Higher order couplings  extend to next nearest neighbours~$U_{i\pm 2}$\,, and beyond, and contain terms of higher order in $r$~and~$\gamma$; typically, $\cos\ell>0.6$\,, even for high orders of coupling~\cite{codevariD}. 

The physical problem of interest is at full coupling $\gamma=1$\,. 
In this case the coupling condition~\eqref{eq:cc} is effectively a Taylor series expansion of the macroscale fields about the centre of the action region $j=\pm rH=\pm(n-a)h$ (for left and right cases), set equal to the average of the microscale fields in the left and right action region~\cite{Roberts11, codevariD}. 
However, centre manifold theory provides full physical support for patch dynamics within some domain about $\gamma=0$ (\(\gamma=0\) is the no coupling case where patches have no influence on each other~\cite{Roberts03, Roberts07, Samaey10}).
Generally we must consider the full range $0\leq \gamma\leq 1$ when providing theoretical support for our patch dynamics scheme; however, as we here avoid formal definitions of $f_{i,\pm n}$~and~$\cos\ell$ we do not delve into a detailed analysis of the limiting behaviour of coupling parameter~$\gamma$.

Section~\ref{sec:micro} determines all eigenvalues and eigenvectors on an arbitrary single patch with patch half-width~$n$, core half-width~$a$, and the coupling conditions~\eqref{eq:cc}.
From these we construct the microscale field solution of~\eqref{eq:ode} on the $i$th~patch with the coupling conditions~\eqref{eq:cc}.  
Importantly, in this solution on the $i$th~patch, the form of the coupling to neighbouring patches, represented by~$f_{i,\pm n}$\,, and the coupling to the $i$th patch, represented by $\cos\ell$, is arbitrary and thus is valid for any form of inter-patch coupling. 
Once the microscale solution on the $i$th~patch is determined, core averaging~\eqref{eq:amp} provides the macroscale solution~$U_i$\,.

In the standard patch dynamics scheme, the coupling conditions~\eqref{eq:cc} are evaluated at each microscale time step as required by the microscale simulator---here the \textsc{ode}~\eqref{eq:ode}.
On a computer with massive parallelisation, and in the scenario where one processor simulates the system~\eqref{eq:ode} over one patch, when applying coupling conditions~\eqref{eq:cc} the coupling data must be transferred between processors each microscale time step. 
As discussed in Section~\ref{sec:intro}, frequent data transfers defeat massive parallelisation.
We propose to limit the transfer of data between processors by only communicating data between patches on mesoscale time-steps~$\delta t$, larger than the microscale time-steps of the simulation but smaller than the macroscale times of interest~$\Delta t$.
Thus, as a first approximation we replace coupling conditions~\eqref{eq:cc} with
\begin{equation} 
\sum_{j=\pm n-2a}^{\pm n}\frac{u_{j+iN}(m\delta t + t)}{2a+1}
=U_i(m\delta t+ t)\cos\ell
+\frac{f_{i,\pm n}(m\delta t)}{2a+1}\,,\label{eq:cc2}
\end{equation}
where $m=0,1,\ldots,M$\,, $0\leq t<\delta t$\,. 
Thus data transfers between processors are required much less frequently: the cost is an error in the simulation which we analyse in Section~\ref{sec:error}.
Figure~\ref{fig:elements} illustrates this new scheme in the case of nearest neighbour coupling.

More sophisticated mesoscale coupling conditions than~\eqref{eq:cc2} are obtained by approximating $f_{i,\pm n}(t)$ in coupling condition~\eqref{eq:cc} with the first~$Q$ terms of its Taylor series expansion about the previous mesoscale time step~$m\delta t$.
Using \(f^q\)~to denote the \(q\)th~derivative of coupling function~\(f\), this generalises the infrequent coupling conditions~\eqref{eq:cc2} to
\begin{equation} 
\sum_{j=\pm n-2a}^{\pm n}\frac{u_{j+iN}(m\delta t + t)}{2a+1}
=U_i(m\delta t+ t)\cos\ell
+\sum_{q=0}^{Q-1} \frac{f_{i,\pm n}^q(m\delta t)t^q}{(2a+1)q!}\,,\label{eq:cc3}
\end{equation}
where the Taylor series expansion is $(Q-1)$th~order accurate.
The error in simulations with such a \(Q\)th~order coupling is also analysed in Section~\ref{sec:error}

Section~\ref{sec:meso} modifies the microscale field solution with coupling conditions~\eqref{eq:cc} obtained by Section~\ref{sec:micro} on the arbitrary $i$th~patch to a solution with mesoscale coupling conditions \eqref{eq:cc2}~or~\eqref{eq:cc3}.
We use these new solutions to systematically explore the errors of various coupling schemes.
Section~\ref{sec:error} analyses the error which arises when coupling conditions~\eqref{eq:cc} are replaced with the mesoscale coupling conditions \eqref{eq:cc2}~or~\eqref{eq:cc3} and its dependence on parameters such as the patch half-width~$n$, core half-width~$a$, ratio~$r$, and Taylor series order of accuracy~$Q-1$.
We consider both the error on the microscale lattice within the $i$th~patch and the error of the macroscale solution~$U_i$ obtained from core averaging~\eqref{eq:amp}.

\section{Microscale sub-patch dynamics}
\label{sec:micro}

To assess the error due to coupling at infrequent mesoscale time steps \eqref{eq:cc2} or~\eqref{eq:cc3},  this section establishes the microscale solution within one arbitrary patch with coupling at microscale time steps~\eqref{eq:cc}.
Section~\ref{sec:meso} modifies the solution with microscale coupling~\eqref{eq:cc} into a solution with mesoscale coupling \eqref{eq:cc2} or~\eqref{eq:cc3}.

To construct the microscale solution within the $i$th~patch we first rewrite equations \eqref{eq:ode}~and~\eqref{eq:cc} as a single matrix equation and then determine all eigenvalues and right and left eigenvectors.
The microscale solution is a linear combination of terms in these eigenvalues and eigenvectors.
This solution is valid for any patch half-width~$n$ and any equal sized action regions and core, $0\leq a<n$\,, with the exception of some special cases where repeated eigenvalues are associated with linearly dependent eigenvectors and the set of all eigenvectors do not form a complete basis on the patch.
These special cases require generalised eigenvectors to provide a full microscale solution~\cite{axler1997}; however, we do not consider generalised eigenvectors here as they further complicate the problem without providing any additional insights.

The case when the microscale is the diffusion \textsc{pde}, will be obtained as the limit of \(n\to\infty\) with \(h\to0\) and finite patch width~$2nh$.

\subsection{Matrix form}
\label{sec:matrix}

In matrix form, the system of \textsc{ode}s~\eqref{eq:ode} within the $i$th~patch and with coupling conditions~\eqref{eq:cc} is
\begin{equation}
B\dot{\vec{u}}(t)=\mathcal{L}\vec{u}(t)+\vec{f}(t) \,,\label{eq:mat}
\end{equation} 
where $(2n+1)$~dimensional vector $\vec{u}=(u_{-n+iN},\ldots,u_{n+iN})$ describes the field~$u_{j+iN}$ at every sub-patch coordinate $j=-n,-n+1,\ldots,n-1,n$\,, and has initial condition $\vec{u}(0)=\vec{u}_0=(u_{-n+iN}(0),\ldots,u_{n+iN}(0))$\,.
The forcing vector $\vec{f}=(f_{i,-n}\,,0,\ldots,0,f_{i,n})$ where $f_{i,\pm n}$ describes the coupling of the $i$th~patch to neighbouring patches, such as the nearest neighbour coupling~\eqref{eq:nn}.
Matrices $B=\operatorname{diag}(0,1,1,\ldots,1,0)$ and $\mathcal{L}$ are $(2n+1)\times(2n+1)$\,.
Rather than use the usual matrix numbering of rows and columns (i.e., $1,2,\ldots,(2n+1)$), we index the rows and columns of $B$~and~$\mathcal{L}$ as $-n,-n+1,\ldots,n-1,n$\,, since these correspond to our patch indices $j=-n,-n+1,\ldots,n-1,n$\,, and similarly for the indexing of components of the vectors $\vec{u}$~and~$\vec{f}$.
For example, $\mathcal{L}_{-n,-n}$~is the element in the first row and first column of~$\mathcal{L}$.

With the exception of the first and last rows, the  nonzero elements of matrix~$\mathcal{L}$ are $\mathcal{L}_{j,j-1}\,, \mathcal{L}_{j,j+1}=1$\,, $\mathcal{L}_{j,j}=-2$ for $j\neq \pm n$\,, which describes the discrete diffusion in \textsc{ode}~\eqref{eq:ode}.
The first and last rows of~$\mathcal{L}$ represent the patch coupling conditions~\eqref{eq:cc} with the macroscale fields  defined by the average over the core in equation~\eqref{eq:amp}.
For an action and core half-width~$a$, $\mathcal{L}_{\pm n,j}=\mathcal{L}^a_{\pm n,j}+\mathcal{L}^c_j$\,, where 
\begin{equation}
\mathcal{L}^a_{\pm n,j}=\begin{cases}
-1, &  n-2a\leq \pm j\leq n\,,\\
0 & \text{otherwise,}
\end{cases}\quad 
\mathcal{L}^c_{j}=\begin{cases}
\cos\ell\,, &  -a\leq j\leq a\,,\\
0 & \text{otherwise.}
\end{cases}\label{eq:c}
\end{equation}
The action and core half-width is only restricted by the size of the patch, $0\leq a<n$\,.
The action regions overlap the core when $a\geq n/3$\,, which means some microscale fields appear both in the average over the core~\eqref{eq:amp} and in the averages over the action regions of the coupling conditions~\eqref{eq:cc}.
Equation~\eqref{eq:c} and the subsequent microscale solution remain valid whether or not the action regions and core overlap. 
For example, when $a=n/3$ the core and action regions overlap at two microscale lattice points~$\pm n/3+iN$ and $\mathcal{L}_{\pm n,\pm n/3}=-1+\cos\ell$\,.

We now solve the generalised eigenproblem and adjoint eigenproblem 
\begin{equation}
(\mathcal{L}-\lambda_k B)\vec{v}_k=\vec{0} \quad \text{and}\quad\vec{z}_k^T(\mathcal{L}-\lambda_k B)=\vec{0}^T,\label{eq:mateqns}
\end{equation}
for right and left eigenvectors $\vec{v}_k$~and~$\vec{z}_k^T$, respectively, and eigenvalues~$\lambda_k$\,. 
We normalise all eigenvectors such that $\vec{z}_k^TB\vec{v}_{k'}=\delta_{kk'}$ for Kronecker delta~$\delta_{kk'}$ and all $k,k'=0,1,\ldots,2n-2$\,. 
The eigenproblem produces two set of, at most, $(2n-1)$~linearly independent right and left eigenvectors indexed $k=0,1,\ldots,2n-2$\,, although the associated $(2n-1)$~eigenvalues may have multiplicity greater than one.
A set of $(2n-1)$~linearly independent right or left eigenvectors forms a complete basis which spans the subspace of the microscale field vector~$\vec{u}$ satisfying coupling conditions~\eqref{eq:cc}.\footnote{The size of the basis which spans the subspace of the microscale field vector~$\vec{u}$ is the number of elements of~$\vec{u}$ (i.e., the number of microscale fields on one patch~$(2n+1)$) minus the number of constraints on the elements of~$\vec{u}$ . There are two constraints, provided by the two coupling conditions~\eqref{eq:cc}, so the size of the basis is~$(2n-1)$.}
Some special cases are degenerate and do not provide the required $(2n-1)$~linearly independent eigenvectors: for simplicity we limit analysis to the generic case of $(2n-1)$~linearly independent eigenvectors.

\subsection{Act and sample on a lattice point}
\label{sec:act}
We first consider the simplest case where there is no averaging over the action region and core, $a=0$\,; instead we just act and sample the microscale field at end points and the mid point of a patch.
The patch coupling conditions~\eqref{eq:cc} only constrain the microscale fields~$u_{j+iN}$ at the patch edges, $j=\pm n$\,, and the macroscale field value in this $i$th~patch~$U_i$ is the microscale value at the centre of the patch, $j=0$\,:  equation~\eqref{eq:amp} reduces to $U_i(t)=u_{iN}(t)$.
The only nonzero elements of the first and last rows of~$\mathcal{L}$ are $\mathcal{L}_{-n,-n}=\mathcal{L}_{n,n}=-1$ and $\mathcal{L}_{\pm n,0}=\cos\ell$\,.
In this case we always obtain $(2n-1)$~linearly independent right or left eigenvectors and $(2n-1)$~distinct eigenvalues, and thus, for this case, the sets of all right or left eigenvectors forms a complete basis (we never need generalised eigenvectors). 

Using Matlab we evaluate the eigenvalues and eigenvectors of the matrix equations~\eqref{eq:mateqns} for several $n$~and~$a$.
From these numerical examples we determine the general analytic forms of the eigenvalues and eigenvectors and confirm by substitution into the matrix equations.
The eigenvalues of the matrix equations~\eqref{eq:mateqns} are
\begin{equation}
\lambda_k=-2\left[1-\cos(\pi l_k/2n)\right]\,,\label{eq:lam0}
\end{equation}
for eigenvector mode indices $k=0,1,\ldots,2(n-1)$ and corresponding wave number
\begin{equation}
l_k=\begin{cases} k+1+(-1)^{k/2}(2\ell/\pi-1)& \text{for even }k,\\
k+1 & \text{for odd }k.\label{eq:lk}
\end{cases}
\end{equation}
The right and left eigenvectors for all sub-patch coordinates~$j$, are, for even~$k$,
\begin{align}
(\vec{v}_k)_j&=(n\sin \ell)^{-1/2} \cos(jl_k\pi/2n)\,,\label{eq:vec0even}\\
(\vec{z}_k)_j&=(n\sin \ell)^{-1/2} \sin[\ell-(-1)^{k/2}|j|l_k\pi/2n]\,,\quad \text{except }(\vec{z}_k)_{\pm n}=(\vec{z}_k)_{\pm (n-1)}\,,\nonumber
\end{align}
(or, equivalently $(\vec{z}_k)_j=(-1)^{k/2}(n\sin \ell)^{-1/2}\sin[(n-|j|)l_k\pi/2n]$ for $j\neq \pm n$)
and for odd~$k$,
\begin{align}
(\vec{v}_k)_j&=n^{-1/2}\sin(jl_k\pi/2n)\,,\nonumber\\
(\vec{z}_k)_j&=n^{-1/2}\sin(jl_k\pi/2n)\,,\quad \text{except }(\vec{z}_k)_{\pm n}=(\vec{z}_k)_{\pm (n-1)}\,\label{eq:vec0odd}.
\end{align}
These eigenvectors are normalized such that $\vec{z}_k^TB\vec{v}_{k'}=\delta_{kk'}$ for all $k,k'=0,1,\ldots,2(n-1)$\,. 


\begin{figure}
\centering
\includegraphics{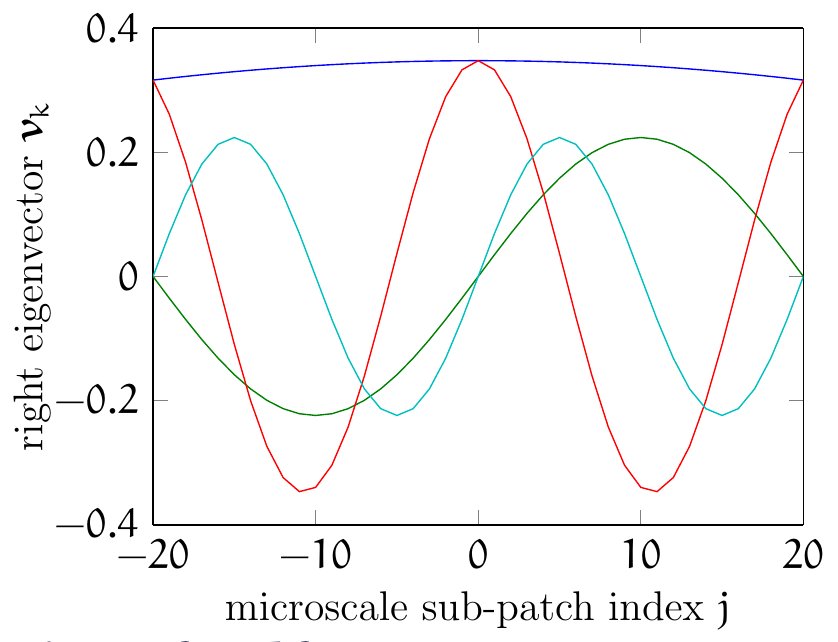}
\includegraphics{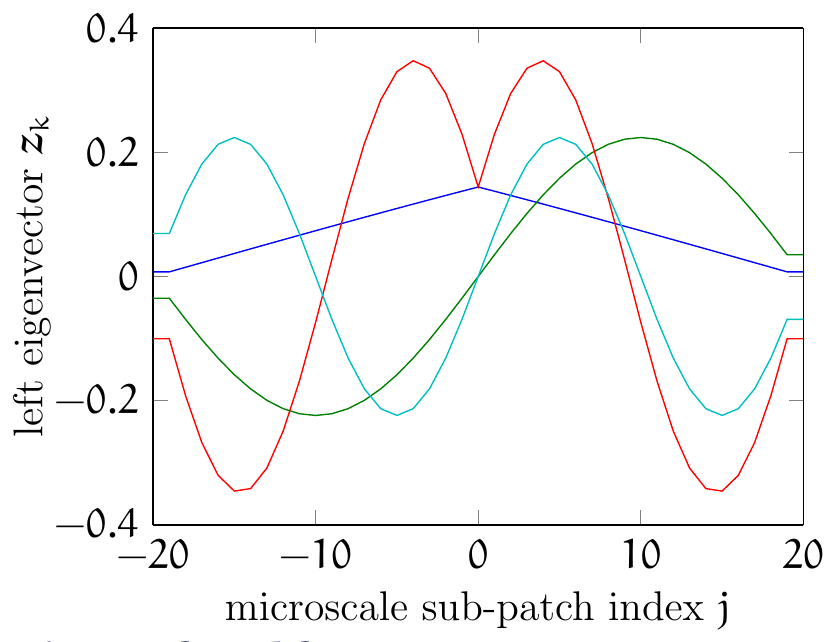}
\caption{Right and left microscale eigenvectors associated with the four lowest magnitude  microscale eigenvalues ($k=0,1,2,3$ for blue, green, red and cyan, respectively) for patch half-width $n=20$\,, action and core half-width $a=0$ and $\cos\ell=0.91$ (so $(n\sin\ell)^{-1/2}=0.35$ and $(\sin\ell/n)^{1/2}=0.14$\,).
}\label{fig:nobuff}
\end{figure}

Figure~\ref{fig:nobuff} plots right and left eigenvectors associated with the four lowest magnitude eigenvalues, for patch half-width $n=20$ and $a=0$\,.
These eigenvectors are only defined at discrete patch lattice points $j=0,\pm 1,\ldots,\pm n$\,, but, for clarity, we plot curves rather than points. 
These plots are typical for any~$n$.
For odd~$k$ the right and left eigenvectors are identical sine functions which are odd about $j=0$ and typical of what one sees in simple one dimensional diffusion~\cite{haberman2012}. 
For even~$k$ the right and left eigenvectors are even functions about $j=0$ but, with their dependence on the parameter~$\ell$, are not typical of eigenvectors of simple one dimensional diffusion.
For all even~$k$ the eigenvectors'  amplitudes are no more than~$(n\sin\ell)^{-1/2}$; however, in the centre of the patch, $j=0$, the right and left eigenvectors behave differently with respect to~$\ell$. 
For even~$k$ at $j=0$ the right eigenvector is $(n\sin\ell)^{-1/2}$ but the left eigenvector is $(\sin\ell/n)^{1/2}$.
So, for small~$\ell$ the right eigenvector is large at $j=0$ but the left eigenvector is small, but as $\ell$~increases the right eigenvector at $j=0$ decreases and the left eigenvector at $j=0$ increases.

Section~\ref{sec:mfcc} uses these, and subsequent eigenvectors, in spectral expansions of the dynamics within the $i$th~patch in order to discover the effects of different mesoscale coupling procedures (Section~\ref{sec:meso}).

\subsection{Action regions and core average over part of a patch}

Let's now consider the case where the core and action regions have half-width $a>0$\,. 
As in Section~\ref{sec:act}, we evaluate numerical examples of eigenvalue and eigenvectors using Matlab and from these we determine the analytic forms, which we confirm by substitution into the matrix equations~\eqref{eq:mateqns}.
The eigenvalues of matrix equations~\eqref{eq:mateqns} are 
\begin{equation}
\lambda_k=\begin{cases}
-2\left\{1-\cos[\pi l_k/2(n-a)]\right\}, &k=0,1,\ldots,2(n-a-1)\,,\\
-2\left\{1-\cos[\pi l_k/(2a+1)]\right\}, &k=2(n-a)-1,\ldots,2(n-1)\,.
\end{cases}
\end{equation}
The $k\leq 2(n-a-1)$ eigenvalues are similar to those for the $a=0$ case, except that the half-width of the patch is now effectively~$(n-a)$ rather than~$n$, thus reducing the buffering of the macroscale solution from $n$~to~$n-a$.
For $k=0,1,\ldots, 2(n-a-1)$ the wavenumbers are the same as the $a=0$ case,
\begin{equation}
l_k=\begin{cases} k+1+(-1)^{k/2}(2\ell/\pi-1)& \text{for even }k,\\
k+1 & \text{for odd }k.\label{eq:lkb}
\end{cases}
\end{equation}

\begin{figure}
\centering
\includegraphics{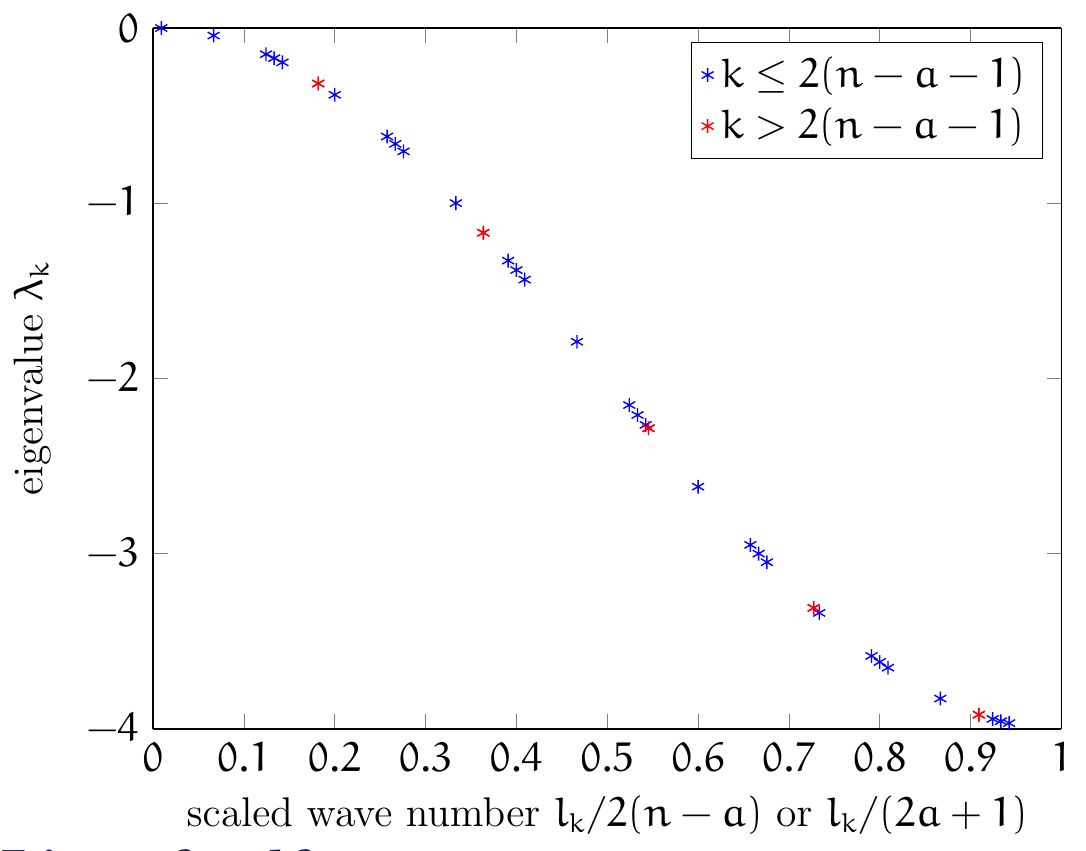}
\caption{Microscale eigenvalues plotted against wavenumbers which are scaled by the reduced width of the patch, $l_k/2(n-a)$, for $k\leq 2(n-a-1)$ (blue) and scale by the width of the core~$l_k/(2a+1)$, for $k> 2(n-a-1)$ (red).
The patch half-width is $n=20$\,, the action and core half-width is $a=5$ and $\cos\ell=0.91$\,.
}\label{fig:wn}
\end{figure}

For $k=m+2(n-a-1)$ where $m=1,2,\ldots, 2a$\,, the wave numbers are 
\begin{equation}
l_{m+2(n-a-1)}=2\lceil m/2 \rceil\,.
\end{equation}
The last $2a$ eigenvalues are $a$~equal pairs. 
If we set $a=0$\,, then the above eigenvalues reduce to those given in equation~\eqref{eq:lam0}.
Figure~\ref{fig:wn} plots scaled wavenumbers against eigenvalues for patch half-width $n=20$ and action and core half-width  $a=5$. 
For this case there are $2(n-a-1)+1=29$ unique eigenvalues and $a$~pairs of equal eigenvalues.

The right eigenvectors for all sub-patch coordinates~$j$ and $k=0,1,\ldots,2(n-a-1)$ are
\begin{equation}
(\vec{v}_k)_j=\begin{cases} [(n-a)\sin\ell]^{-1/2}\cos[jl_k\pi/2(n-a)]& \text{for even }k,\\
(n-a)^{-1/2}\sin[jl_k\pi/2(n-a)]& \text{for odd }k.
\end{cases}
\end{equation}
If we set $a=0$\,, then the above eigenvectors reduce to the right eigenvectors given in equations \eqref{eq:vec0even}~and~\eqref{eq:vec0odd}.
For $k=m+2(n-a-1)$ where $m=1,2,\ldots,2a$\,, the right eigenvectors for all sub-patch coordinates~$j$ are
\begin{equation}
(\vec{v}_k)_j=\begin{cases} (2a+1)^{-1/2}\cos[jl_k\pi/(2a+1)]& \text{for even }k,\\
(2a+1)^{-1/2}\sin[jl_k\pi/(2a+1)]& \text{for odd }k.
\end{cases}
\end{equation}

Despite there being $a$~pairs of equal eigenvalues for $k>2(n-a-1)$, the associated right eigenvectors are linearly independent since for even~$k$ the eigenvectors are even in~$j$ and for odd~$k$ the eigenvectors are odd in~$j$.\footnote{It is possible, for some particular $n$~and~$a$, for there to be odd $k_1\leq 2(n-a-1)$ and odd $k_2>2(n-a-1)$ such that $\lambda_{k_1}=\lambda_{k_2}$ which implies $\vec{v}_{k_1}=\vec{v}_{k_2}$\,. 
For example, such a case occurs for $n=4$ and $a=1$ when $k_1=3$ and $k_2=5$\,.
Such repeated eigenvectors mean that we do not have a complete basis. 
To form a complete basis we should replace repeated eigenvectors with generalised eigenvectors~\cite[e.g.]{axler1997},  but we  avoid this complication by avoiding those $n$~and~$a$ which result in repeated eigenvectors.}

The left eigenvectors are considerably more complex than the right eigenvectors. 
To define the left eigenvectors we first define some new functions.
For $k=0,1,\ldots,2(n-a-1)$\,,
\begin{align}
(\vec{w}_k)_j={}&(-1)^{\lceil(k-1)/2\rceil}\{2\sin[(2a+1)l_k\pi/4(n-a)]\}^{-1}\nonumber\\
{}&\times\left[\cos\left(\frac{l_k\pi}{4(n-a)}\right)-\cos\left(\frac{(2j+1)l_k\pi}{4(n-a)}\right)\right]\nonumber\\
{}&\times \begin{cases}
[(n-a)\sin\ell]^{-1/2} & \text{for even }k,\\
(n-a)^{-1/2} & \text{for odd }k,
\end{cases}\label{eq:w1}
\end{align}
and for $k=m+2(n-a-1)$ where $m=1,2,\ldots,2a$\,,
\begin{align}
(\vec{w}_k)_j={}&(-1)^{\lceil(m-2)/2\rceil}\left[\cos\left(\frac{l_k\pi}{2(2a+1)}\right)-\cos\left(\frac{(2j+1)l_k\pi}{2(2a+1)}\right)\right]\label{eq:w2}\\
{}&\times\begin{cases}
(2a+1)^{1/2}\{\sin[(n-a)l_k\pi/(2a+1)]\}^{-1}, &\text{for odd }m,\\
(2a+1)^{1/2}\{\cos[(n-a)l_k\pi/(2a+1)]-\cos\ell\}^{-1}, & \text{for even }m. 
\end{cases}\nonumber
\end{align}
We also define, for $k=0,1,\ldots,2(n-a-1)$ and $|j|=0,1,\ldots,n-a$\,,
\begin{equation}
(\vec{z}'_k)_j=\begin{cases} 
[(n-a)\sin\ell]^{-1/2}\sin[\ell-(-1)^{k/2}|j|l_k\pi/2(n-a)]\,, & \text{even }k,\\
(n-a)^{-1/2}\sin[jl_k\pi/2(n-a)]\,, & \text{odd }k,
\end{cases}
\end{equation}
and for all other possible values of $k$~and~$j$ we set $(\vec{z}'_k)_j=0$\,.

The left eigenvectors for sub-patch coordinate $|j|\leq n-1$ are, for even~$k$,
\begin{equation}
(\vec{z}_k)_j=(\vec{z}'_k)_j+\begin{cases}
-2\cos\ell(w_k)_{a-|j|}\,, & 0\leq|j|\leq a \,, \\
0, & b\leq |j|\leq n-2a \,, \\
(w_k)_{|j|-n+2a}\,, & n-2a\leq|j|\leq n-a \,, \\
(w_k)_{n-|j|}\,, & n-a\leq |j|\leq n-1 \,,
\end{cases}
\end{equation}
and for odd~$k$,
\begin{equation}
(\vec{z}_k)_j=(\vec{z}'_k)_j+\begin{cases}
0, & 0\leq |j|\leq n-2a \,, \\
\sgn(j)(w_k)_{|j|-n+2a}\,, & n-2a\leq|j|\leq n-a \,, \\
\sgn(j)(w_k)_{n-|j|}\,, & n-a\leq |j|\leq n-1 \,.
\end{cases}
\end{equation}
In all cases, $(\vec{z}_k)_{\pm n}=(\vec{z}_k)_{\pm (n-1)}$\,.
These~$\vec{z}_k$ reduce to the left eigenvectors given in equations \eqref{eq:vec0even} and~\eqref{eq:vec0odd} for the special case \(a=0\)\,.

Equation~\eqref{eq:w1} is undefined when $k=k_1\leq 2(n-a-1)$ is odd and $(2a+1)l_{k_1}\pi/4(n-a)=s\pi$ for integer $s$, or equivalently, when $l_{k_1}\pi/2(n-a)=l_{k_2}\pi/(2a+1)$ for $k_2=m+2(n-a-1)$ and $m=2s-1$\,.
From the limits on~$k_1$ it can be shown that $s=1,2,\ldots,a$ so $m=1,3,\ldots,(2a-1)$\,.
Therefore, equation~\eqref{eq:w1} is undefined if for odd $k=k_1<2(n-a-1)$  there exists odd $k_2>2(n-a-1)$ such that $\lambda_{k_1}=\lambda_{k_2}$\,; a scenario which we identify as a degenerate case which requires generalised eigenvalues and is avoided.
Similarly, equation~\eqref{eq:w2} is undefined when $k=k_2>2(n-a-1)$ is odd and $(n-a)l_{k_2}\pi/(2a+1)=s\pi$ for integer~$s$, or equivalently, when $l_{k_2}\pi/(2a+1)=l_{k_1}\pi/2(n-a)$ for $k_1=2s-1$\,.
From the limits on~$k_2$ it can be shown that $s=1,2,\ldots,(n-a-1)$ so $k_1=1,3,\ldots,2(n-a-1)-1$\,.
Therefore, equation~\eqref{eq:w2} is undefined if for odd $k=k_2>2(n-a-1)$ there exists odd $k_1\leq 2(n-a-1)$ such that $\lambda_{k_2}=\lambda_{k_1}$\,, which is, again, a degenerate case to be avoided.

\begin{figure}
\centering
\includegraphics{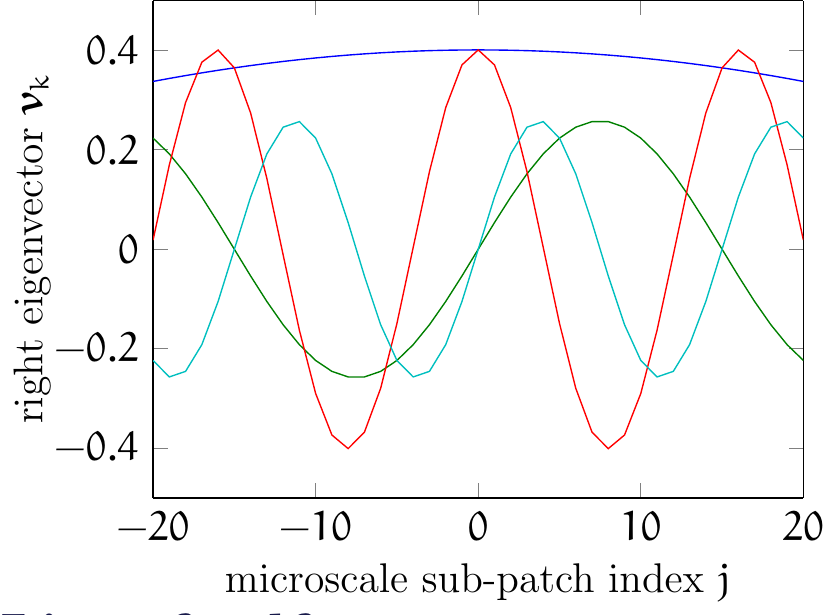}
\includegraphics{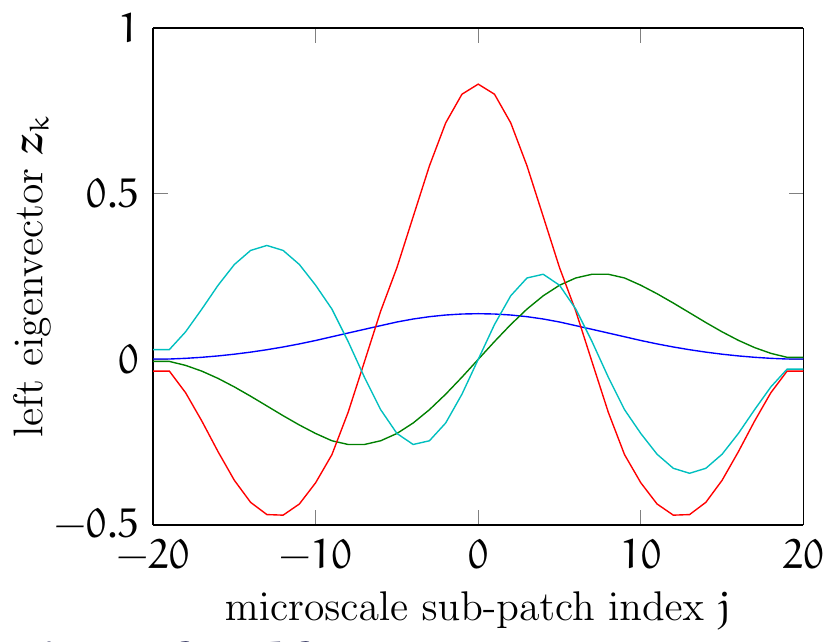}
\caption{Right and left microscale eigenvectors associated with the four lowest magnitude eigenvalues with $k\leq 2(n-a-1)$ ($k=0,1,2,3$ for blue, green, red and cyan, respectively) with the same parameters as Figure~\ref{fig:wn}, patch half-width $n=20$\,, action and core half-width $a=5$ and $\cos\ell=0.91$\,.}\label{fig:patch}
\end{figure}

\begin{figure}
\centering
\includegraphics{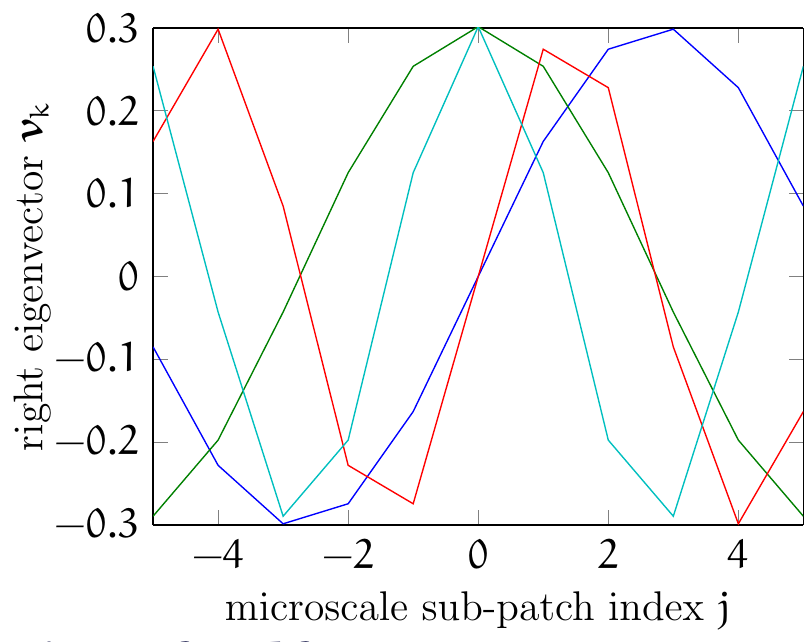}
\includegraphics{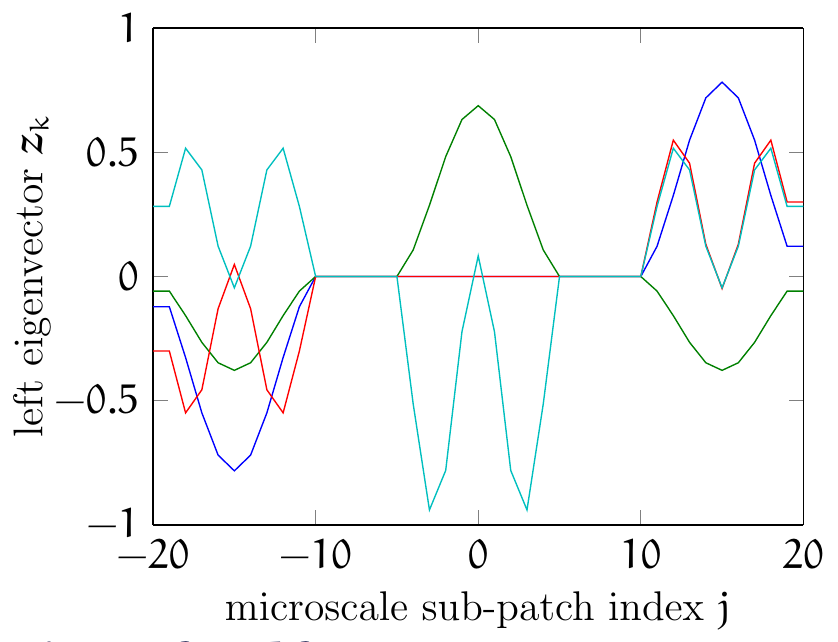}
\caption{Right and left microscale eigenvectors associated with the four lowest magnitude eigenvalues with $k>2(n-a-1)$ ($k=m+2(n-a-1)$ with $m=1,2,3,4$ for blue, green, red and cyan, respectively) with the same parameters as Figure~\ref{fig:wn}, patch half-width $n=20$\,, action and core half-width $a=5$ and $\cos\ell=0.91$\,.
The right eigenvectors are only plotted over the patch core, $|j|\leq a$\,, as this shows at least one complete period. 
The left eigenvectors are plotted over the entire patch, $|j|\leq n$\,.}\label{fig:buff}
\end{figure}

Figure~\ref{fig:patch} plots right and left eigenvectors associated with the four lowest magnitude eigenvalues with $k\leq 2(n-a-1)$\,, for patch half-width $n=20$ and $a=5$\,, and Figure~\ref{fig:buff} plots right and left eigenvectors associated with the four lowest magnitude eigenvalues with $k> 2(n-a-1)$\,, with the same patch geometry.
These plots are typical of any~$n$ and $a\neq 0$\,, provided we have a complete basis of eigenvectors across sub-patch coordinates $|j|=0,1,\ldots,n-1$ (so do not require generalised eigenvectors).
The shape of the right eigenvectors are not remarkably different from those for $a=0$ shown in Figure~\ref{fig:nobuff}, with the main point of difference being the frequency of the sinusoidal eigenvectors which are increased by the effective reduction of patch size from~$n$ to~$n-a$.
In contrast, the left eigenvectors for $a>0$ are significantly different to those for $a=0$\,.
For nonzero action regions and core, $\vec{z}_k$~with odd~$k\leq 2(n-a-1)$ appears smooth about $j=0$\,, unlike~$\vec{z}_k$ for zero action regions and core. 
For $k> 2(n-a-1)$ the left eigenvectors with odd~$k$ are only nonzero inside the action regions whereas for even~$k$ the left eigenvectors are only nonzero inside the action regions and the core. 

\subsection{Microscale field for continuous coupling}
\label{sec:mfcc}

Given the homogeneous microscale solutions of the previous subsection, we now explore a spectral representation of the solution within a patch with forcing by general coupling with neighbouring patches.
The microscale field solution of \textsc{ode}~\eqref{eq:mat} within the $i$th~patch, which has continuous and instantaneous coupling with nearby patches defined by coupling conditions~\eqref{eq:cc}, is
\begin{align}
\vec{u}(t)&=\sum_{k=0}^{2(n-1)}e^{\lambda_kt}\vec{v}_k\vec{z}_k^TB\vec{u}(0)+\vec{f}(t)+\sum_{k=0}^{2(n-1)}\vec{v}_k\vec{z}_k^T\int_0^t\vec{f}(t')e^{\lambda_k(t-t')}\,dt'\nonumber\\
&=T(t)B\vec{u}_0+\vec{f}(t)+T(t)\star \vec{f}(t)\,,\label{eq:gen}
\end{align}
where $\vec{u}=(u_{-n+iN},\ldots,u_{n+iN})$ describes all microscale fields~$u_{j+iN}$ with $|j|\leq  n$ within the $i$th~patch, and the convolution in the last term is defined as
\begin{equation}
T(t)\star \vec{f}(t)=\sum_{k=0}^{2(n-1)}T_k\int_0^t\vec{f}(t')e^{\lambda_k(t-t')}\,dt'\,,
\end{equation}
with $(2n+1)\times(2n+1)$ matrices $T_k=\vec{v}_k\vec{z}_k^T$ and state transition matrix
\begin{equation}
T(t)=\sum_{k=0}^{2(n-1)}e^{\lambda_k t}T_k\,.
\end{equation}
This solution is confirmed via direct substitution into \textsc{ode}~\eqref{eq:mat} and with the identity, proved in Appendix~\ref{app:identity}, that
\begin{equation}
T(0)=\sum_{k=0}^{2(n-1)}T_k=\sum_{k=0}^{2(n-1)}\vec{v}_k\vec{z}^T_k=B+A\,,\label{eq:T0}
\end{equation}
where all elements of the $(2n+1)\times(2n+1)$ matrix~$A$ are zero, except 
\begin{equation}
A_{j,\pm n}=\mathcal{L}_{j,\pm n}\quad\text{and} \quad A_{\pm n,j}=\mathcal{L}_{\pm n,j}\quad \text{for }|j|=0,1,\ldots,n\,.\label{eq:A}
\end{equation}
Thus, on substituting solution~\eqref{eq:gen} into \textsc{ode}~\eqref{eq:mat} we use
\begin{equation}
B\sum_{k=0}^{2(n-1)}T_k\vec{f}(t)=BA\vec{f}(t)=(0,f_{i,-n},0,\ldots,0,f_{i,n},0)=\vec{f}(t)+\mathcal{L}\vec{f}(t)\,.
\end{equation}
Identity~\eqref{eq:T0} is also useful in confirming solution~\eqref{eq:gen} has the correct initial condition. 

The general forced solution~\eqref{eq:gen} acts as the reference for reporting errors due to infrequent mesoscale coupling.

\section{Infrequent mesotime coupling}
\label{sec:meso}

This section approximates the exact analytic solution~\eqref{eq:gen} of \textsc{ode}~\eqref{eq:ode} within one patch with coupling conditions~\eqref{eq:cc} by replacing the continuous time coupling vector $\vec{f}(t)$ with a mesoscale coupling vector.
With mesoscale coupling we only evaluate the coupling vector~$\vec{f}(t)$ at fixed time intervals $t=m\delta t$ where the meso-time interval \(\delta t\) is much larger than the time-step of the microscale computation.
We assume we know the first $Q$~terms of the Taylor series in time of the coupling~$\vec{f}(t)$ exactly at times $t=m\delta t$ for $m=0,1,\ldots,M$\,, but at no other points in time.
If we only know the values of~\(\vec f\), and not their time derivatives, then parameter \(Q=1\)\,; but our analysis applies for general~\(Q\).
Further, we assume the computation in a patch up to time $t=m\delta t$  only depends on~$\vec{f}(t')$ evaluated at $t'<m\delta t$\,. 
Thus we know the $(Q-1)$th~order Taylor series of $\vec{f}(m\delta t+t)$ about $t=0$ for $0\leq t<\delta t$ and we apply coupling conditions~\eqref{eq:cc3}.
We first consider the initial time step with $m=0$ where the coupling~$\vec{f}(t)$ and its derivatives are only known at $t=0$\,.
By homogeneity in time,  the analysis extends straightforwardly to a general number of meso-time steps.

\subsection{Initial time step}

Consider the solution~\eqref{eq:gen} at $t=\delta t$ multiplied by matrix~$B$,
\begin{equation}
B\vec{u}(\delta t)
=B\sum_{k=0}^{2(n-1)}\mu_kT_kB\vec{u}_0+\sum_{k=0}^{2(n-1)}BT_k\vec{J}_{k}\,,\label{eq:Bu}
\end{equation}
where $\mu_k=e^{\lambda_k\delta t}$ and
\begin{equation}
\vec{J}_{k}=\int_0^{\delta t}\vec{f}(t')e^{\lambda_k(\delta t-t')}\,dt'.\label{eq:ibp0}
\end{equation}
The multiplication by~$B$ in the above solution removes the fields at the endpoints of the patch,~$u_{\pm n}$, which can always be evaluated from the coupling conditions~\eqref{eq:cc}, \eqref{eq:cc2}~or~\eqref{eq:cc3}, if required, provided all other field components are known. 
The aim is to approximate~$\vec{J}_k$ in~\eqref{eq:Bu} using the $(Q-1)$th~order Taylor series expansion of $\vec{f}(t')$, as used in coupling conditions~\eqref{eq:cc3}, and also to provide the error of~$\vec{J}_k$ due to this Taylor series approximation. 

Integrate~\eqref{eq:ibp0} by parts,
\begin{equation}
\vec{J}_{k}=\frac{1}{\lambda_k}\vec{f}(0)(\mu_k-1)+
\frac{1}{\lambda_k}\int_0^{\delta t}\dot{\vec{f}}(t')(e^{\lambda_k (\delta t-t')}-1)dt'\,,\label{eq:ibp1}
\end{equation}
where a constant of integration is chosen so that only~$\vec{f}(0)$, not~$\vec{f}(\delta t)$, appears in the integrated part.
The first term in equation~\eqref{eq:ibp1} is exactly what is obtained by substituting the zeroth order Taylor series of $\vec{f}(t')$ about $t'=0$ into integral~\eqref{eq:ibp0} (i.e., replace~$\vec{f}(t')$ with~$\vec{f}(0)$) and thus the integral part of~\eqref{eq:ibp1} is the error of the zeroth order Taylor series approximation.
Then, 
\begin{equation}
B\vec{u}(\delta t)
=\sum_{k=0}^{2(n-1)}BT_k\left[\mu_kB\vec{u}_0+\vec{f}(0)(\mu_k-1)/\lambda_k\right]+B\vec{R}\,,\label{eq:u1}
\end{equation}
with remainder vector $\vec{R}=(R_{-n},R_{-n+1},\ldots,R_n)$ where, for $|j|=0,1,\ldots,(n-1)$\,, 
\begin{align}
R_j={}&\sum_{k=0}^{2(n-1)}\frac{1}{\lambda_k}\int_0^{\delta t}\left[BT_k\dot{\vec{f}}(t')\right]_j(e^{\lambda_k (\delta t-t')}-1)dt'\nonumber\\
={}&\sum_{k=0}^{2(n-1)}\frac{(v_{k})_j}{\lambda_k}\int_0^{\delta t}\left[(z_{k})_{-n}\dot{f}_{-n}(t')+(z_{k})_{n}\dot{f}_{n}(t')\right](e^{\lambda_k (\delta t-t')}-1)dt'\,.\label{eq:Rj}
\end{align}
Since the microscale fields in the action regions are constrained by the coupling conditions~\eqref{eq:cc2}, the remainders at $j=\pm n$ are dependent on the other remainders in the action regions,
\begin{equation}
R_{\pm n}=-\sum_{j=\pm n-2a}^{\pm (n-1)}R_j\,.\label{eq:Redge}
\end{equation}
In the mesoscale coupling scheme with coupling conditions~\eqref{eq:cc2} we use
\begin{equation}
B\vec{u}(\delta t)
\sim \sum_{k=0}^{2(n-1)}BT_k\left[\mu_kB\vec{u}_0+\vec{f}(0)(\mu_k-1)/\lambda_k\right]\,,\label{eq:approxu1}
\end{equation}
with the error known to be precisely~$B\vec{R}$\,.

The solution~\eqref{eq:ibp1} for~$\vec{J}_k$ is for one integration by parts which is equivalent to approximating $\vec{f}(t')$ by a zeroth order Taylor series ($Q=1$), as in coupling conditions~\eqref{eq:cc2}. 
Generalising to $Q$~integrations by parts,
\begin{align}
\vec{J}_{k}={}&\sum_{q=1}^Q\vec{f}^{q-1}(0)\lambda_k^{-q}\left(\mu_k-\sum_{p=0}^{q-1}\delta t^p\lambda_k^p/p!\right)\nonumber\\
{}&+\frac{1}{\lambda_k^{Q}}\int_0^{\delta t}\vec{f}^{Q}(t')\left[e^{\lambda_k(\delta t-t')}-\sum_{p=0}^{Q-1}\frac{(-1)^p}{p!}\lambda_k^{p}(t'-\delta t)^{p}\right]dt'\,,\label{eq:ibpQ}
\end{align}
where we choose constants of integration so that only~$\vec{f}(0)$ and its derivatives, not~$\vec{f}(\delta t)$ and its derivatives, appear in the integrated part.
In this case
\begin{equation}
B\vec{u}(\delta t)
=\sum_{k=0}^{2(n-1)}BT_k\left[\mu_kB\vec{u}_0+\sum_{q=1}^Q\frac{\vec{f}^{q}(0)}{\lambda_k^{q}}\left(\mu_k-\sum_{p=0}^{q}\frac{\delta t^p\lambda_k^p}{p!}\right)\right]+B\vec{R}\,,\label{eq:sol1Q}
\end{equation}
with components of the remainder vector
\begin{align}
R_j={}&\sum_{k=0}^{2(n-1)}\frac{1}{\lambda_k^{Q}}\int_0^{\delta t}[BT_k\vec{f}^{Q}(t')]_j\left[e^{\lambda_k(\delta t-t')}-\sum_{p=0}^{Q-1}\frac{(-1)^p}{p!}\lambda_k^{p}(t'-\delta t)^{p}\right]dt'\nonumber\\
={}&\sum_{k=0}^{2(n-1)}\frac{1}{\lambda_k^{Q}}\int_0^{\delta t}(\vec{v}_k)_j[(\vec{z}_k)_{-n}f_{-n}^{Q}(t')+(\vec{z}_k)_{n}f_n^{Q}(t')]\sum_{p=Q}^{\infty}\frac{(-1)^p}{p!}\lambda_k^{p}(t'-\delta t)^{p}\,dt'\nonumber\\
={}&\sum_{p=0}^{\infty}\left[\sum_{k=0}^{2(n-1)}\frac{\lambda_k^{p}}{(Q+p)!}(\vec{v}_k)_j(\vec{z}_k)_{-n}\right]\int_0^{\delta t}f_{-n}^{Q}(t')(\delta t-t')^{Q+p}\,dt'\nonumber\\
{}&+\sum_{p=0}^{\infty}\left[\sum_{k=0}^{2(n-1)}\frac{\lambda_k^{p}}{(Q+p)!}(\vec{v}_k)_j(\vec{z}_k)_{n}\right]\int_0^{\delta t}f_{n}^{Q}(t')(\delta t-t')^{Q+p}\,dt'\,,\label{eq:cj}
\end{align}
for $|j|=0,1,\ldots,(n-1)$, and where $R_{\pm n}$ satisfy~\eqref{eq:Redge}.

On substituting the $(Q-1)$th~order Taylor series of~$\vec{f}(t')$ into integral~\eqref{eq:ibp0}, a general term in the expansion, for $q=0,\ldots,(Q-1)$\,, is
\begin{equation}
\frac{1}{q!}\int_0^{\delta t}\vec{f}^q(0)t'^qe^{\lambda_k(\delta t-t')}\,dt'=\vec{f}^q(0)\lambda_k^{-(q+1)}\left(\mu_k-\sum_{p=0}^q\delta t^p\lambda_k^p/p!\right),
\end{equation}
which describes all terms in the first line of the expansion of~$\vec{J}_k$ in~\eqref{eq:ibpQ}.
Thus, this first line of~\eqref{eq:ibpQ} is the approximation of~$\vec{J}_k$ with mesoscale coupling conditions~\eqref{eq:cc3} obtained from the $(Q-1)$th order Taylor series expansion of $\vec{f}(t)$ and the second line of~\eqref{eq:ibpQ} is the error of~$\vec{J}_k$ due to the Taylor series expansion.
In the mesoscale coupling scheme, when we know $(Q-1)>0$~derivatives of~$\vec{f}(t)$ at $t=0$ we improve on the mesoscale coupling conditions~\eqref{eq:cc2} and approximate field solution~\eqref{eq:approxu1} by using coupling conditions~\eqref{eq:cc3} which provides the approximate field solutions
\begin{equation}
B\vec{u}(\delta t)
\sim\sum_{k=0}^{2(n-1)}BT_k\left[\mu_kB\vec{u}_0+\sum_{q=1}^Q\vec{f}^{q-1}(0)\lambda_k^{-q}\left(\mu_k-\sum_{p=0}^{q-1}\delta t^p\lambda_k^p/p!\right)\right]\label{eq:approxuQ}
\end{equation}
with the error known to be precisely~$B\vec{R}$\,.

\subsection{Multiple mesotime steps}
\label{sec:Msteps}

After $M$ mesoscale time steps of size~$\delta t$, solution~\eqref{eq:gen} gives
\begin{equation}
B\vec{u}(M\delta t)=B\sum_{k=0}^{2(n-1)}\mu_k^MT_kB\vec{u}(0)
+\sum_{k=0}^{2(n-1)}BT_k\sum_{m=0}^{M-1}\mu_k^{M-m-1}\vec{J}_{km}
\end{equation}
where the integral from time zero to~$M\delta t$ is converted into $M$~integrals from zero to~$\delta t$,
\begin{equation}
\vec{J}_{km}=\int_0^{\delta t}\vec{f}(t'+m\delta t)e^{\lambda_k(\delta t-t')}\,dt'\,.\label{eq:ibpm}
\end{equation}
The integral~$\vec{J}_{km}$ is a generalised version of the integral~$\vec{J}_k$ in equation~\eqref{eq:ibp0}.
After $Q$~integrations by parts, $\vec{J}_{km}$ is similar to~$\vec{J}_k$ in equation~\eqref{eq:ibpQ}, but with $\vec{f}(0)$ and~$\vec{f}(t')$ replaced with $\vec{f}(m\delta t)$ and~$\vec{f}(t'+m\delta t)$\,, respectively.
Thus we obtain
\begin{equation}
B\vec{u}(M\delta t)=\sum_{k=0}^{2(n-1)}BT_k\left[\mu_k^MB\vec{u}(0)
+\sum_{q=1}^Q\frac{\vec{f}_M^{q-1}(0)}{\lambda_k^{q}}\left(\mu_k-\sum_{p=0}^{q-1}\frac{\delta t^p\lambda_k^p}{p!}\right)\right]+B\vec{R}\label{eq:solMQ}
\end{equation}
where 
\begin{equation}
\vec{f}_M^q(t')=\sum_{m=0}^{M-1}\mu_k^{M-m-1}\vec{f}^q(t'+m\delta t)\,,\label{eq:fM}
\end{equation}
the components of the remainder vector for $|j|=0,1,\ldots,(n-1)$ are
\begin{equation}
R_j=\sum_{k=0}^{2(n-1)}\frac{1}{\lambda^{Q}_k}
\int_0^{\delta t}[BT_k\vec{f}_M^Q(t')]_j\left[e^{\lambda_k(\delta t-t')}-\sum_{p=0}^{Q-1}\frac{(-1)^p}{p!}\lambda^p_k(t'-\delta t)^p\right]dt'\,,\label{eq:cjM}
\end{equation}
and $R_{\pm n}$ satisfy~\eqref{eq:Redge}.

The solution over $M$~mesoscale time steps~\eqref{eq:solMQ} is the same as the solution over one mesoscale time step~\eqref{eq:sol1Q}, but with derivatives of the coupling vector~$\vec{f}(0)$ replaced by derivatives of the $M$~time steps coupling vector~$\vec{f}_M(0)$.
The only difference between the remainder vector~$\vec{R}$ over $M$~mesoscale time steps and over one mesoscale times step, is that for~$M$, equation~\eqref{eq:cjM} integrates over~$\vec{f}_M^Q(t')$ and  for one time step, equation~\eqref{eq:cj} integrates over~$\vec{f}^Q(t')$\,.
In the following error analysis we only consider a single mesoscale time step, but these results are readily adaptable to multiple mesoscale time steps.

\section{Error analysis}
\label{sec:error}

The remainder vector~$\vec{R}$~\eqref{eq:cj},  or~\eqref{eq:cjM} for multiple mesoscale time steps, tells us the extent to which coupling errors penetrate into the core of the $i$th~patch from the action regions.
Here we are only concerned with errors arising from evaluating the patch coupling~$\vec{f}(t)$ at mesoscale time intervals in the mesoscale coupling conditions \eqref{eq:cc2} or~\eqref{eq:cc3}, not with errors which are due to, for example, the interpolation between patches. 
Ideally, the remainder~$R_j$ should be minimal within the patch core (i.e., for~$|j|\leq a$) since it is the microscale fields~$u_j$ within the patch core which determine the macroscale field~$U_i(t)$.
Errors near the patch edges (i.e., $j\approx \pm n$) need not be small. 
Section~\ref{sec:rj}  considers how~$R_j$ varies across the patch, $|j|=0,1,\ldots,n$\,, and with varying the number of terms~$Q$ in the Taylor series expansion in the coupling conditions~\eqref{eq:cc3}.
Section~\ref{sec:U} then looks at the upper bound of the error of the macroscale field~$U_i(\delta t)$ due to mesoscale coupling.
This analysis considers one mesoscale time step as the base to predict multiple steps. 

\subsection{Error penetration in patch}
\label{sec:rj}

Defining $\bar{f}^Q(\delta t)=\max_{t'\in[0,\delta t]}\big[|f_{i,-n}^Q(t')|,|f_{i,n}^Q(t')|\big]$ and rearranging equation~\eqref{eq:cj} gives the scaled upper error bound for the remainder when $|j|\leq (n-1)$\,,
\begin{equation}
\frac{|R_j|}{\bar{f}^Q(\delta t)}\leq {}\frac{\delta t^{Q+1}}{(Q+1)!}\sum_{k=0}^{2(n-1)}\left|(\vec{v}_k)_j\left[(\vec{z}_k)_{-n}+\vec{z}_k)_{n}\right]{}_1F_1(1;Q+2;\lambda_k\delta t)\right|=R_{j\max}\,,\label{eq:scaler}
\end{equation}
where the sum over~$p$ is rewritten as a confluent hypergeometric function~\cite{Grad2014}
\begin{equation}
\sum_{p=0}^{\infty}\frac{\lambda_k^{p}}{(Q+p)!}\int_0^{\delta t}(\delta t-t')^{Q+p}dt'
= \frac{\delta t^{Q+1}}{(Q+1)!}{}\,_1F_1(1;Q+2;\lambda_k\delta t)\,.
\end{equation}
For $j=\pm n$\,, from equation~\eqref{eq:Redge},
\begin{equation}
R_{\pm n\max}=\sum_{j=\pm n-2a}^{\pm (n-1)}R_{j\max}\,.\label{eq:Redgemax}
\end{equation}

\begin{figure}
\centering
\includegraphics{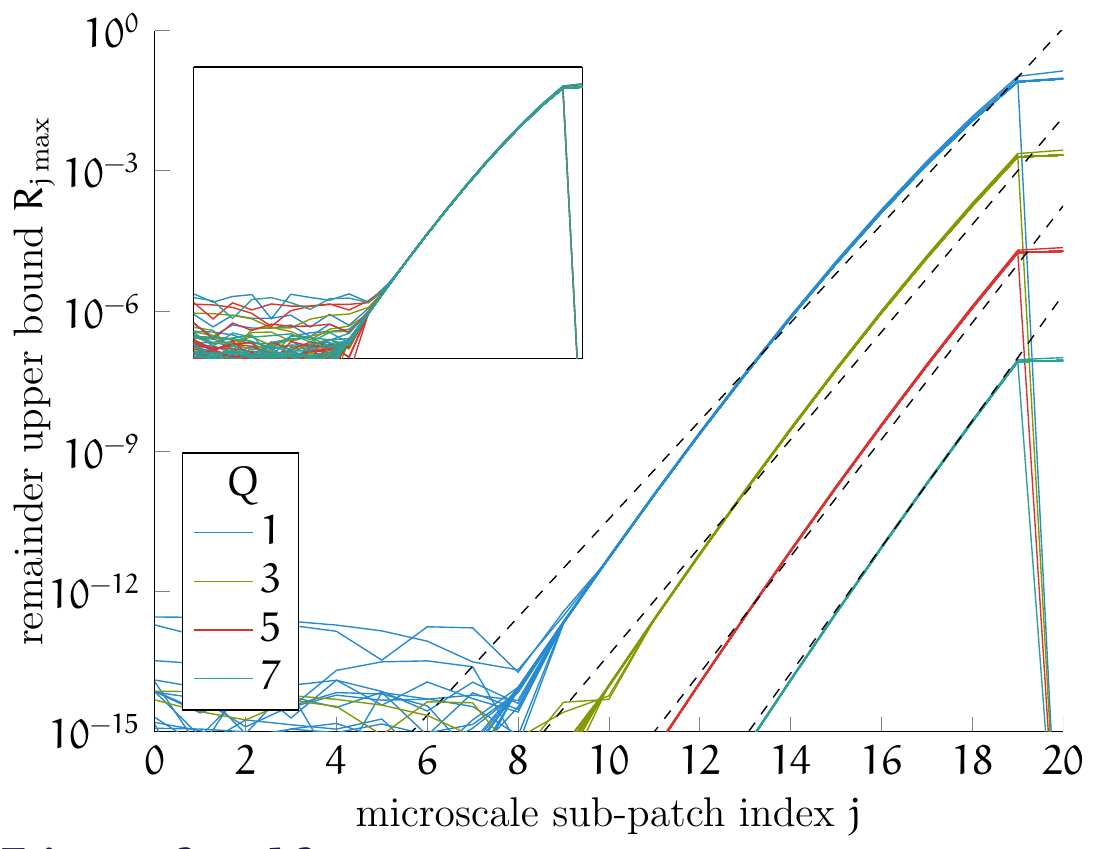}
\caption{The solid lines describe the upper bound of components of the remainder vector~$R_{j\max}$ in a patch with patch half-width $n=20$\,,  $\cos\ell=0.91$\,, mesoscale time~$\delta t=0.5$  and $Q=1,3,5,7$\,.
Lines with the same colour have the same~$Q$ but all possible core half-widths $a=0,\ldots,19$\,; these different core half-widths are only distinguishable in the centre of the patch $j<10$\,, and on the very edge of the patch $j=20$ for the $a=0$ case.
The dashed lines approximate the coloured curves with a simple power function.
Numerical error dominates at low~$R_{j\max}$\,.
Inset:~the remainder vector~$R_{j\max}$ on the same scale as the main plot with $n=20$\,, $\delta t=0.5$\,, $Q=1$\,, and $\cos\ell=0.65, 0.75, 0.85, 0.95$\,. 
Lines with the same colour have the same~$\cos\ell$ but all possible core half-widths $a=0,\ldots,19$ and are only distinguishable below $R_{j\max}\sim 10^{-13}$ where numerical error dominates.}\label{fig:penetrate20}
\end{figure}

Figure~\ref{fig:penetrate20} plots the upper bound of the remainder~$R_{j\max}$ for a range of~$Q$ and $\delta t=0.5$\,.
Since the remainder is symmetric about $j=0$\,, only the $j\geq 0$ components are shown.
Figure~\ref{fig:penetrate20} shows that the remainder is very small in the core of the patch, as required.
For $\delta t\sim 0.5$ the upper bound is approximated by
\begin{equation}
R_{j\max}\sim (2\delta t)^{Q+1} 10^{-Q-(1+0.025Q/\delta t)(n-1-j)}\label{eq:approx}
\end{equation} 
until numerical error dominates around $R_{j\max}\sim10^{-13-Q}$\,.
With constant patch half-width~$n$ and $|j|<n$, the core half-width~$a$ has little affect on the remainder~$R_{j\max}$ and, until numerical error is dominant, there is nothing to distinguish the remainders with different~$a$, provided the mesoscale time~$\delta t$ and number of terms in the Taylor series expansion~$Q$ are fixed.
At the patch edge $|j|=n$ the remainder upper bounds for core half-widths $a>0$ are distinctly different from the  $a=0$ case.
For $a>0$\,, in~\eqref{eq:Redgemax} the $R_{\pm(n-1)\max}$~term dominates the sum and $R_{\pm n\max}\sim R_{\pm(n-1)\max}$\,.
For $a=0$ the sum in equation~\eqref{eq:Redgemax} vanishes and thus $R_{\pm n\max}=0$\,.

The inset of Figure~\ref{fig:penetrate20} shows that varying~$\cos\ell$, which is a function of the ratio $r=(n-a)/N$ and the coupling strength~$\gamma$, does not affect the upper bound of the remainder.  
Furthermore, as implied by the approximation~\eqref{eq:approx}, the upper bound of the remainder is dependent on the distance from the patch edge, that is $n-j$\,, rather than the patch half-width~$n$ alone, and~$R_{j\max}$ decreases exponentially with $n-j$.
Therefore, larger~$n$ does not in general produce smaller~$R_{j\max}$; however, it does produce wider patches which have smaller~$R_{j\max}$ within the centre core of these patches. 
In macroscale modelling the most important error measurement is the error in the macroscale field, and the macroscale field is obtained by averaging over microscale fields in the patch core.
Therefore, for the system considered here the errors in the macroscale field will generally be minimised for larger patch half-widths~$n$ and smaller core half-widths~$a$.

We expect similar results for other microscale systems in the same universality class as~\eqref{eq:ode}.

\subsection{Macroscale solution}
\label{sec:U}

The macroscale field value obtained from the $i$th~patch,~$U_i(\delta t)$, is the average of the microscale field~$u_{j+iN}(\delta t)$ in the patch core, as shown in equation~\eqref{eq:amp}.
The error of this macroscale value is the average of the remainder vector components~$R_j$ in the patch core:
\begin{align}
E=\frac{1}{2a+1}\left|\sum_{j=-a}^aR_j\right|
={}&\left|\frac{1}{2a+1}\sum_{p=0}^{\infty}\left[\sum_{k=0,\text{even}}^{2(n-1)}\frac{\lambda_k^{p}}{(Q+p)!}\sum_{j=-a}^a(\vec{v}_k)_j(\vec{z}_k)_{n}\right]\right.\nonumber\\
{}&\left.\times\int_0^{\delta t}[f_{-n}^Q(t')+f_{n}^Q(t')](\delta t-t')^{Q+p}dt'\right|,
\end{align}
where only the even eigenvalues and eigenvectors  contribute. 

Similarly to the scaled upper bound of the remainder in~\eqref{eq:scaler}, the scaled upper bound of the error is
\begin{equation}
\frac{(2a+1)E}{\bar{f}^Q(\delta t)}\leq \frac{\delta t^{Q+1}}{(Q+1)!}\sum_{k=0,\text{even}}^{2(n-1)}\sum_{j=-a}^a\left|(\vec{v}_k)_j(\vec{z}_k)_{n}{}\,_1F_1(1;Q+2;\lambda_k\delta t)\right|=E_{\max}\,.
\end{equation} 
In this error bound we scale by~$(2a+1)$ because $\bar{f}^Q$~is of the order $(2a+1)$ due to scaling of the coupling conditions~\eqref{eq:cc3} and, for example,~\eqref{eq:nn}.

\begin{figure}
\centering
\includegraphics{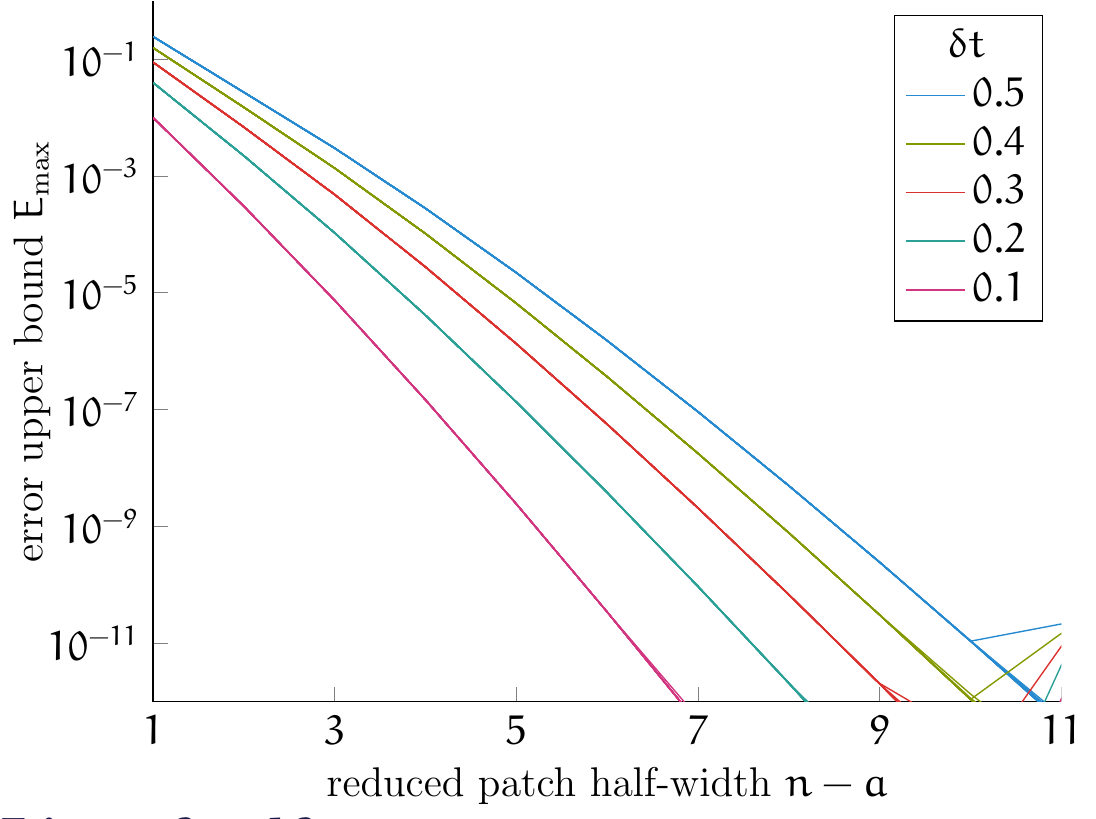}
\caption{The upper bound of the error~$E_{\max}$ for $Q=1$ and $\cos\ell=0.91$\,, a range of  mesoscale time steps~$\delta t$ and over several reduced patch half-widths~$n-a$\,. The calculated range of the patch and core half-widths was for all $4\leq n\leq 20$ and $0\leq a<n$\,, respectively, but since the error~$E_{\max}$ is extremely small for $n-a\geq 11$ these cases are not plotted.
Curves with different~$n$ and~$a$ are identical when they have the same reduced patch half-width~$n-a$ and mesoscale time step~$\delta t$ (and are thus plotted in the same colour), until numerical error dominates the calculation below $E_{\max}\sim 10^{-11}$.}\label{fig:er}
\end{figure}

Figure~\ref{fig:er} plots the scaled error bound~$E_{\max}$ for $Q=1$ (communicating only function values~\(f_{\pm n}(0)\) and no derivatives as in coupling conditions~\eqref{eq:cc2}) and a range of mesoscale time-steps~$\delta t$.
This figure shows that the upper bound of error~$E_{\max}$ is a function of reduced patch half-width~$n-a$, rather than individual values of $n$~and~$a$, and decreases as $n-a$ increases.
So, the error is smaller for large~$n$ and small~$a$, as predicted at the end of Section~\ref{sec:rj}.
Section~\ref{sec:rj} also showed that the remainder upper bound~$R_{j\max}$ is independent of~$\cos\ell$, and thus independent of ratio $r=(n-a)/N$\,, and we conclude that the error upper bound~$E_{\max}$ should then also be independent of $\cos\ell$~and~$r$.
This is a surprising result as it implies that the error is independent of the ratio between the microscale and macroscale lattice spacings $N=H/h$\,. 
We interpret this to mean that it is most important to make the reduced patch half-width~$n-a$ (or buffer width) only large enough to capture a significant portion of the microscale dynamics, and then macroscale modelling may proceed to any scale. 
The interpretation of `significant portion' will depend on the nature of the microscale model.
Of course, here we are exploring an error bound scaled by~$\bar{f}^Q(\delta t)$;  some~$r$ dependence will be present in the absolute error when, for example, we use~$f_{i,\pm}(t)$ of the form~\eqref{eq:nn}. 
 
Even for only one term in the Taylor series expansion of the coupling condition~\eqref{eq:cc2}, $Q=1$\,, the upper bound of error~$E_{\max}$ plotted in Figure~\ref{fig:er} decreases significantly with increasing~$n-a$ or decreasing~$\delta t$.
Based on Figure~\ref{fig:penetrate20} we expect larger values of~$Q$ to produce smaller errors for macroscale field~$U_i(t)$, but in a simulation, increasing~$Q$ might not be practical.
This is because a given~$Q$ requires the first~$Q-1$ temporal derivatives of the coupling vector~$\vec{f}$, and this information may not be practically available. 
Whether it is more practical to increase~$n-a$ or decrease~$\delta t$ in order to reduce errors is dependent on the problem being considered and the architecture of the parallel computer running the simulation. 
For example, in cases where there is periodic microscale detail, choosing $n-a$~to be multiples of the period gives more accuracy~\cite{codevariD} and so it might be best to choose an optimal~$n-a$ and then consider reducing~$\delta t$.
In general, while increasing $n-a$ will increase processing times, decreasing~$\delta t$ will increase the amount of data communication. 
In the case of a supercomputer with a large number of processors, maximising~$n-a$ will take advantage of the processing power while choosing~$\delta t$ large enough will avoid limitations associated with interprocessor communication.

\section{Two dimensional Ginzburg--Landau numerical simulation}
\label{sec:numerics}

To demonstrate that mesoscale coupling is effective in dimensions higher than one and in more complicated problems than the simple diffusion of \textsc{ode}~\eqref{eq:ode}, we use patch dynamics with mesoscale coupling to simulate the complex Ginzburg--Landau equation~\cite{Aranson2002} on a discrete two dimensional square microscale spatial lattice with lattice spacing~$h$:
\begin{align}
\dot{u}_{j_x,j_y}(t)={}&{}(1+i\alpha)[u_{j_x+1,j_y}(t)+u_{j_x-1,j_y}(t)-2u_{j_x,j_y}(t)]\nonumber\\
{}&{}+(1+i\alpha)[u_{j_x,j_y+1}(t)+u_{j_x,j_y-1}(t)-2u_{j_x,j_y}(t)]\nonumber\\
{}&{}+u_{j_x,j_y}(t)-(1+i\beta) u_{j_x,j_y}(t)|u_{j_x,j_y}(t)|^2,\label{eq:gl}
\end{align}
for constant, real $\alpha$~and~$\beta$.
The complex Ginzburg--Landau equation contains two dimensional diffusion terms with an additional nonlinear part.
We set $\alpha=1$ and $\beta=2$\,; this choice of parameters should, after some time, produce plane wave solutions. 
The initial conditions specified were sinusoidal with amplitude~$0.5$ plus a real, normally distributed random number with mean~$0$ and standard deviation~$0.8$ at each microscale lattice point.

The construction of a two dimensional macroscale lattice with lattice spacing~$H$ and two dimensional patches is very similar to the one dimensional case discussed in Section~\ref{sec:patchsetup}.
As in the one dimensional case, we assume that $N=H/h$ is an integer. 
We construct the square macroscale lattice with general lattice point~$(X_{i_x},Y_{i_y})=(i_xH,i_yH)$ for integers $i_{x,y}$\,.
Centred about each macroscale lattice point~$(X_{i_x},Y_{i_y})$ we construct the $(i_x,i_y)$~square patch of width~$2nh$ which does not touch or overlap neighbouring patches. 
Again, like the one dimensional case, integer~$n$ is the patch half-width and to ensure the patches do not overlap $n/N<1/2$\,.
The microscale fields on the $(i_x,i_y)$~patch are $u_{j_x+i_xN,j_y+i_yN}$ for microscale sub-patch index $j_{x,y}=0,\pm1,\ldots,\pm n$\,.
In the patch dynamics numerical simulation the Ginzburg--Landau equation~\eqref{eq:gl} is solved for all microscale values within a patch, excluding those on the patch edges, $j_{x,y}=0,\pm1,\ldots,\pm (n-1)$\,.

We could define the macroscale fields in terms of some average over microscale fields in the centre or core of the patch (i.e., average over fields in a square core with core half-width~$a$), as we did in the one dimensional case~\eqref{eq:amp}.
Similarly, we could define patch coupling conditions in terms of some average over microscale values on the patch edge, as in equations~\eqref{eq:cc}, \eqref{eq:cc2}~or~\eqref{eq:cc3}. 
However, Section~\ref{sec:error} showed that mesoscale errors are a function of the reduced patch half-width~$n-a$ rather than patch half-width~$n$ so here we simply choose the two dimensional analogue of $a=0$ where we act and sample the microscale field at patch edge points and patch mid-points to define patch coupling conditions and macroscale fields, respectively. 
Thus, the macroscale field is the microscale field in the centre of the square patch~\cite{Roberts11}, 
\begin{equation}
U_{i_x,i_y}(t)=u_{i_xN,i_yN}(t)\,,
\end{equation}
and the patch coupling conditions with mesoscale coupling~$\delta t$ constrain each microscale field on four patch edges,
\begin{align}
u_{\pm n+i_xN,j_y}(m\delta t+t)&=U_{i_x,i_y}(m\delta t+t)\cos\ell_{\pm n,j_y}+f_{(i_x,i_y),(\pm n,j_y)}(m\delta t)\,, \nonumber\\
u_{j_x,\pm n+i_yN}(m\delta t+t)&=U_{i_x,i_y}(m\delta t+t)\cos\ell_{j_x,\pm n}+f_{(i_x,i_y),(j_x,\pm n)}(m\delta t)\,, 
\end{align}
for nonnegative integer~$m$ and $j_{x,y}=0,\pm1,\ldots,\pm (n-1)$\,.
These coupling conditions are analogous to the $Q=1$ one dimensional coupling conditions~\eqref{eq:cc2}. 
Both $\cos\ell_{j_x,j_y}$ and $f_{(i_x,i_y),(j_x,j_y)}$ are functions of the patch coupling strength~$\gamma$ and the ratio~$r=nh/H$ and describe interpolation from the centre of the $(i_x,i_y)$~patch and across the un-simulated space between the $(i_x,i_y)$~patch and adjacent patches~\cite{Roberts11}.
Coupling conditions are not required for microscale fields on the patch corners where both $j_x=\pm n$ and $j_y=\pm n$ since these microscale fields are not required when solving the Ginzburg--Landau equation~\eqref{eq:gl} for the microscale fields within the patches but not on the patch edges.

For the numerical simulations presented here we use nearest neighbour coupling derived from Taylor series expansions about the microscale points on the patch edge~\cite{Roberts11}: 
\begin{align}
\cos\ell_{j_x,j_y}={}&{}1-(r_{j_x}^2-r_{j_y}^2)\gamma\,,\nonumber\\
f_{(i_x,i_y),(j_x,j_y)}(t)={}&{}\tfrac12 r_{j_x}\gamma[(r_{j_x}\pm 1)U_{i_x+1,i_y}(t)+(r_{j_x}\mp 1)U_{i_x-1,i_y}(t)]\nonumber\\
{}&{}+\tfrac12 r_{j_y}\gamma[(r_{j_y}\pm 1)U_{i_x,i_y+1}(t)+(r_{j_y}\mp 1)U_{i_x,i_y-1}(t)]\,,
\end{align}
for $r_{j_x}=j_xh/H$ and $r_{j_y}=j_yh/H$ where $j_{x,y}=0,\pm1,\ldots,\pm n$ and $r_{\pm n}=\pm r$\,.
We define a square periodic domain of width~$20$ with microscale lattice spacing $h=0.25$\,.
We set macroscale lattice spacing $H=5$\,, which fits $4\times 4$~patches across the domain and gives $N=H/h=20$\,.
We choose a patch half-width $n=6$ so that $r=0.3$ and choose mesoscale coupling~$\delta t=0.2$\,.

\begin{figure}
\centering
\includegraphics{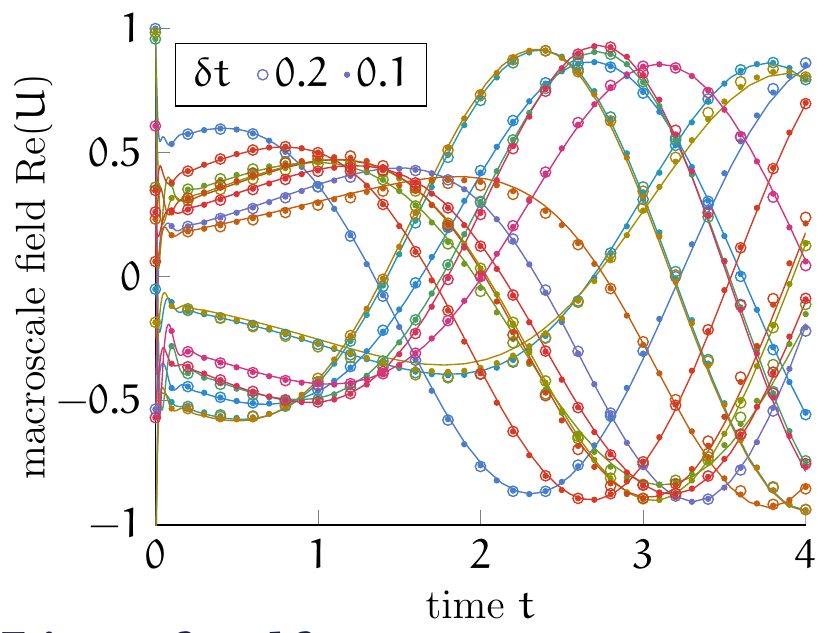}
\includegraphics{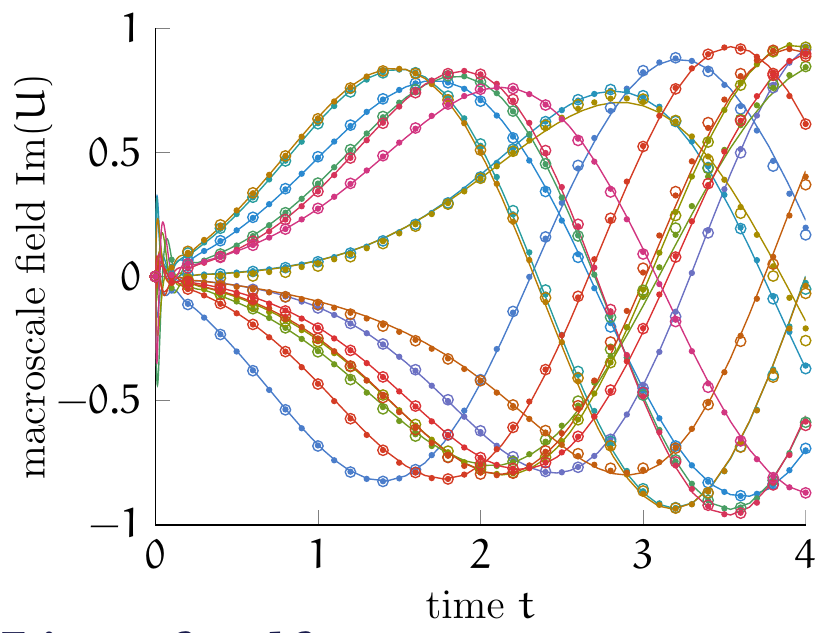}
\caption{Real and imaginary parts of macroscale fields~$U_{i_x,i_y}$ from all $16$~patches with patch half-width $n=6$\,.  Compare continuous time coupling (solid lines) and mesoscale coupling plotted at the mesoscale time steps $\delta t=0.2$ (open circles) and  $\delta t=0.1$ (dots).
}\label{fig:realimag}
\end{figure}

Figure~\ref{fig:3d} compares the real components of the microscale fields within patches for continuous time coupling and with mesoscale coupling at two times, $t=0.04$ and $t=0.4$\,; there is little difference between these two cases. 
The large random component in the initial condition decays rapidly and is substantially reduced even at small time $t=0.04$\,.
By $t=0.4$ the simulation is smooth.
Figure~\ref{fig:realimag} compares the real and imaginary part of the macroscale fields~$U_{i_x,i_y}$ obtained with continuous time coupling and mesoscale couplings $\delta t=0.2, 0.1$\,.
The $\delta t=0.2$ case produces a reasonable result, but some macroscale fields are slightly inaccurate. 
Better results are obtained when the mesoscale coupling time-step is reduced to $\delta t=0.1$\,.

\section{Conclusion}
\label{sec:conclude}

We analysed a patch dynamics macroscale modelling scheme which is adapted for massive parallelisation and exascale computing by limiting the transfer of data between processors, assuming that one processor only calculates the dynamics of one patch. 
Rather than transferring information concerning coupling conditions at every microscale time step, we propose limiting the data transfer to mesoscale time-steps~$\delta t$.
This method does not require any preprocessing and makes minimal use of stored data. 
Limiting the transfer and size of data addresses several major hurdles facing the development of exascale computing, specifically, slow data transfer speeds compared to processor speeds, the energy cost of data transfer and memory size limitations~\cite{doe2010, doe2013, emwg2014}. 
If one processor was to evaluate the dynamics of a small number of adjacent patches, then coupling data should be updated at microscale times for the patches evaluated on that one processor, but coupling data transfers at mesoscale times should be maintained for patches evaluated on different processors.

Section~\ref{sec:error} found that errors arising from mesoscale coupling are controlled by the reduced patch half-width~$n-a$ and the mesoscale coupling time~$\delta t$, with larger~$n-a$ and smaller~$\delta t$ producing smaller errors. 
However, these error predictions only compare patch dynamics with continuous time coupling to patch dynamics with mesoscale coupling; they do not consider errors inherent in the patch dynamics scheme due extrapolation across large regions of un-simulated space, and these additional errors are also functions of~$n-a$.
Previous work analysed the error of patch dynamics with continuous time coupling relative to the known exact solution and found that larger~$n-a$ are not always better when the model has rough microscale detail~\cite{codevariD}.
If the microscale detail has some periodicity, then optimal solutions require that~$n-a$ is chosen such that the periodicity exactly divides~$n-a$.
Thus, adjusting the reduced patch half-width~$n-a$ for improvements in mesoscale coupling patch dynamics is not practical if~$n-a$ is already constrained by the symmetry of the microscale model. 
In these cases, to reduce errors for mesoscale coupling it is advisable to fix~$n-a$ at the optimal solution determined from the symmetry of the microscale model and only reduce~$\delta t$.

Section~\ref{sec:intro} briefly discusses the need for fault management and resiliency in exascale computing~\cite{doe2010, doe2013, emwg2014}.
While the proposed patch dynamics scheme makes no allowances for failure, we expect that some modifications to the scheme will address this issue.  
Future work may develop coupling conditions for patch dynamics which are dependent on \emph{variable} mesoscale times steps~$\delta t$, so that some delay in data transfer is accounted for in the algorithm. 
In addition, this variable time step may be permitted to extend to infinity, thus accounting for cases where the data never arrives, due to, say, complete failure of a processor. 

The mathematical analysis presented here only considers a simple one dimensional microscale diffusion model, but the scheme is readily modifiable to patches in two or more spatial dimensions and more complex models~\cite{Roberts11}.
The numerical simulation in Section~\ref{sec:numerics} show that mesoscale temporal coupling is effective for the complex two dimensional nonlinear Ginzburg--Landau model, and reducing the mesoscale coupling time~$\delta t$  produces more accurate results, in agreement with the one dimensional analysis.   
More work need to be done to develop the full implementation of patch dynamics for mesoscale coupling, that is, with patches in both space and time~\cite{Kevrekidis09} illustrated in Figure~\ref{fig:patchesSpaceTime}.

\appendix
\section{Proof of identity~\eqref{eq:T0}}
\label{app:identity}
Here we prove the identity~\eqref{eq:T0} which defines the state transition matrix~$T(t)$ at time $t=0$\,.
The proof holds for~$\mathcal{L}$ defined in Section~\ref{sec:matrix}, provided the matrix equations~\eqref{eq:mateqns} produce $k=0,1,\ldots,2(n-1)$ linearly independent right and left eigenvectors, $\vec{v}_k$ and~$\vec{z}^T_k$, which satisfy the normalisation condition $\vec{v}_kB\vec{z}_{k'}^T=\delta_{kk'}$ for all $k,k'=0,1,\ldots,2(n-1)$\,.
The set of all~$\vec{v}_k$ or all~$\vec{z}_k$ form a complete basis which spans the subspace of the microscale field vector on one patch~$\vec{u}$ satisfying the coupling conditions~\eqref{eq:cc}, \eqref{eq:cc2} or~\eqref{eq:cc3}.

Define some vector in the subspace of the microscale field vector~$\vec{u}$ in terms of the right eigenvectors
\begin{equation}
\vec{s}=B\sum_{k=0}^{2(n-1)} c_k\vec{v}_k\,,
\end{equation}
for arbitrary coefficients~$c_k$\,.
The first and last components of this vector are $\vec{s}_{\pm n}=0$\,, but, since the right eigenvectors form a complete basis for sub-patch coordinates~$j=-n+1,-n+2,\ldots,n-1$, all other components of~$\vec{s}$ are completely general.
Using the normalisation condition of the eigenvectors,
\begin{equation}
BT(0)\vec{s}=B\sum_{k,k'=0}^{2(n-1)}\vec{v}_k\vec{z}^T_kBc_{k'}\vec{v}_{k'}=B\sum_{k=0}^{2(n-1)}c_k\vec{v}_k=\vec{s}\,,
\end{equation}
and since the form of~$\vec{s}$ is arbitrary for all but its first and last components we conclude that $BT(0)=B+M$
where the matrix~$M$ is such that $M\vec{s}=\vec{0}$\,.
Since~$\vec{s}$ is a general vector, except for $\vec{s}_{\pm n}=0$\,, the only nonzero elements of~$M$ are in its first and last columns. 
In addition, since all elements in the first and last rows of~$BT(0)$ must be zero, we also know that all elements in the first and last rows of~$M$ must be zero. 
We define some matrix~$A$ which satisfies $M=BA$ and, given the form of~$M$, conclude that the only nonzero elements of~$A$ are in the first and last rows and the first and last columns. 
Therefore,
\begin{equation}
T(0)=B+A\,.
\end{equation}
We now show that the matrix~$A$ is of the form given in equation~\eqref{eq:A}.

From equation~\eqref{eq:mateqns} we obtain $(\vec{z}^T_k\mathcal{L})_{\pm n}=0$\,.
Thus, for $j=-n+1,-n+2,\ldots,n-1$\,,
\begin{align}
[T(0)\mathcal{L}]_{j,\pm n}=0&=\sum_{k=0}^{2(n-1)}[\vec{v}_k\vec{z}^T_k\mathcal{L}]_{j,\pm n}=\sum_{k=0}^{2(n-1)}(\vec{v}_k)_j(\vec{z}^T_k\mathcal{L})_{\pm n}\nonumber\\
&=\mathcal{L}_{j,\pm n}+A_{j,-n}\mathcal{L}_{-n,\pm n}+A_{j,n}\mathcal{L}_{n,\pm n}\,. \label{eq:T0L}
\end{align}
Therefore, for $j\neq \pm n$\,,  using the definition of~$\mathcal{L}$ in Section~\ref{sec:matrix},
\begin{equation}
A_{j,\pm n}=\delta_{j,\pm (n-1)}=\mathcal{L}_{j,\pm n}\,.\label{eq:Acolumn}
\end{equation}
Similarly, from equation~\eqref{eq:mateqns} we obtain $(\mathcal{L}\vec{v}_k)_{\pm n}=0$\,, so, for, $j=-n,-n+1,\ldots,n$\,,
\begin{align}
[\mathcal{L}T(0)]_{\pm n,j}=0&=\sum_{k=0}^{2(n-1)}[\mathcal{L}\vec{v}_k\vec{z}^T_k]_{\pm n,j}=\sum_{k=0}^{2(n-1)}(\mathcal{L}\vec{v}_k)_{\pm n}(\vec{z}^T_k)_j\nonumber\\
&=\mathcal{L}_{\pm n,j}+\mathcal{L}_{\pm n,-n}A_{-n,j}+\mathcal{L}_{\pm n,n}A_{n,j}\,. \label{eq:LT0}
\end{align}
For $j\neq \pm n$ and using equation~\eqref{eq:c},
\begin{equation}
A_{\pm n,j}=\mathcal{L}^1_{\pm n,j}+\mathcal{L}^2_j=\mathcal{L}_{\pm n,j}\,.\label{eq:Arow}
\end{equation}
Finally, since $[\mathcal{L}T(0)]_{\pm n,n}=[\mathcal{L}T(0)]_{\pm n,-n}=0$\,,
\begin{align}
\sum_{l=-n}^n \mathcal{L}_{\pm n,l}A_{l,n}&=\mathcal{L}_{\pm n, n}A_{ n, n}+\mathcal{L}_{\pm n,(n-1)}A_{(n-1),n}=0\,,\nonumber\\
\sum_{l=-n}^n \mathcal{L}_{\pm n,l}A_{l,-n}&=\mathcal{L}_{\pm n,-n}A_{-n,-n}+\mathcal{L}_{\pm n,-(n-1)}A_{-(n-1),-n}=0\,,
\end{align}
so that
\begin{equation}
A_{\pm n,\pm n}=-1=\mathcal{L}_{\pm n,\pm n}\quad \text{and}\quad A_{\pm n,\mp n}=0=\mathcal{L}_{\pm n,\mp n}\,.\label{eq:Acorners}
\end{equation}
Therefore, with equations~\eqref{eq:Acolumn}, \eqref{eq:Arow} and~\eqref{eq:Acorners} we have shown that the matrix~$A$ is of the form given in equation~\eqref{eq:A}.

\bibliography{mesoscalecoup}
    
\end{document}